\documentclass[12pt]{article}
\usepackage{amssymb}
\usepackage{latexsym}
\newenvironment{conjalt}[1]{\vspace{1ex}
\par\noindent{\bf Conjecture} #1 \em}{\vspace{1ex}\par}
\newenvironment{property}[1]{\vspace{1ex}
\par\noindent{\bf Property} #1 \em}{\vspace{1ex}\par}
\newtheorem{claim}{Claim}[section]

\newtheorem{cor}{Corollary}[section]

\newtheorem{lemma}{Lemma}[section]
\newtheorem{prop}{Proposition}[section]

\newtheorem{thm}{Theorem}[section]

\newcounter{example}[section]
\newcounter{rem}[section]

\newcommand{\barbbQ}{\ensuremath{\bar{\bbQ}}}
\newcommand{\bbB}{\ensuremath{\mathbb{B}}}
\newcommand{\bbC}{\ensuremath{\mathbb{C}}}
\newcommand{\bbCp}{\ensuremath{\mathbb{C}_p}}

\newcommand{\bbN}{\ensuremath{\mathbb{N}}}

\newcommand{\bbQ}{\ensuremath{\mathbb{Q}}}

\newcommand{\bbR}{\ensuremath{\mathbb{R}}}
\newcommand{\bbZ}{\ensuremath{\mathbb{Z}}}
\newcommand{\bbZp}{\ensuremath{\mathbb{Z}_p}}

\newcommand{\beql}[1]{\begin{equation}\label{#1}}
\newcommand{\bPf}{\noindent \textsc{Proof\ }}

\newcommand{\cA}{\ensuremath{\mathcal{A}}}

\newcommand{\cD}{\ensuremath{\mathcal{D}}}
\newcommand{\cE}{\ensuremath{\mathcal{E}}}

\newcommand{\cf}{\emph{cf.}}
\newcommand{\cG}{\ensuremath{\mathcal{G}}}

\newcommand{\cK}{\ensuremath{\mathcal{K}}}
\newcommand{\Cl}{{\rm Cl}}
\newcommand{\cL}{\ensuremath{\mathcal{L}}}
\newcommand{\cM}{\ensuremath{\mathcal{M}}}

\newcommand{\cO}{\ensuremath{\mathcal{O}}}

\newcommand{\cR}{\ensuremath{\mathcal{R}}}

\newcommand{\cT}{\ensuremath{\mathcal{T}}}
\newcommand{\cU}{\ensuremath{\mathcal{U}}}

\newcommand{\displaymapdef}[5]
{\[
\begin{array}{rcrcl}
 #1 &:& #2 &\longrightarrow& #3 \\
    & &    &                    \\
    & & #4 &\longmapsto    & #5
\end{array}
\]}

\newcommand{\eeq}{\end{equation}}
\newcommand{\eg}{\emph{e.g.}}
\newcommand{\ePf}{\hspace*{\fill}~$\Box$\vertsp\par}
\newcommand{\eps}{\ensuremath{\varepsilon}}
\newcommand{\etc}{\emph{etc.}}
\newcommand{\example}{\refstepcounter{example}
\noindent{\sc Example \theexample}}
\newcommand{\fa}{\ensuremath{\mathfrak{a}}}

\newcommand{\fc}{\ensuremath{\mathfrak{c}}}
\newcommand{\ff}{\ensuremath{\mathfrak{f}}}
\newcommand{\fg}{\ensuremath{\mathfrak{g}}}
\newcommand{\fh}{\ensuremath{\mathfrak{h}}}
\newcommand{\fm}{\ensuremath{\mathfrak{m}}}

\newcommand{\fp}{\ensuremath{\mathfrak{p}}}
\newcommand{\fP}{\ensuremath{\mathfrak{P}}}
\newcommand{\fq}{\ensuremath{\mathfrak{q}}}

\newcommand{\fw}{\ensuremath{\mathfrak{w}}}

\newcommand{\fz}{\ensuremath{\mathfrak{z}}}
\newcommand{\Gal}{{\rm Gal}}

\newcommand{\half}{\frac{1}{2}}

\newcommand{\ie}{\emph{i.e.}}
\newcommand{\inv}{^{-1}}
\newcommand{\ndiv}{\nmid}
\newcommand{\nin}{\not\in}
\newcommand{\ord}{{\rm ord}}

\newcommand{\refeq}[1]{~(\ref{eq:#1})}
\newcommand{\rem}{\refstepcounter{rem}\noindent{\sc Remark \therem}}

\newcommand{\resp}{\emph{resp.}}

\renewcommand{\theexample}{\thesection.\arabic{example}}
\renewcommand{\therem}{\thesection.\arabic{rem}}

\newcommand{\vertsp}{\vspace{1ex}}

\newcommand{\bareps}{\ensuremath{\bar{\eps}}}
\newcommand{\Clmk}{\ensuremath{\Cl_\fm(k)}}
\newcommand{\Deltad}{\ensuremath{\Delta^\dag}}
\newcommand{\Gammad}{\ensuremath{{\Gamma^\dag}}}
\newcommand{\Gammadd}{\ensuremath{{\Gamma^\ddag}}}
\newcommand{\Ginf}{\ensuremath{{G_{\infty}}}}
\newcommand{\Ginfd}{\ensuremath{{G_\infty^\dag}}}
\newcommand{\Ginfdd}{\ensuremath{{G_\infty^\ddag}}}
\newcommand{\GQbarQ}{\ensuremath{\Gal(\barbbQ/\bbQ)}}
\newcommand{\hchi}{\ensuremath{{\hat{\chi}}}}
\newcommand{\Id}{\ensuremath{{\rm Id}}}
\newcommand{\IwaZG}{\ensuremath{{\rm Iwa(\bbZ\times\Gammadd)}}}
\newcommand{\Kinf}{{\ensuremath{K_{\infty}}}}
\newcommand{\kinf}{{\ensuremath{k_{\infty}}}}
\renewcommand{\L}[1]{\ensuremath{L^{(#1)}}}

\newcommand{\PhigT}{\ensuremath{\Phi_{\fg,T}}}

\newcommand{\PhiKkT}{\ensuremath{\Phi_{K/k,T}}}
\newcommand{\PhiKkTp}{\ensuremath{\Phi_{K/k,T,p}}}
\newcommand{\PhimT}{\ensuremath{\Phi_{\fm,T}}}
\newcommand{\PhimTp}{\ensuremath{\Phi_{\fm,T,p}}}
\newcommand{\PhijmTp}{\ensuremath{\Phi^{(j)}_{\fm,T,p}}}
\newcommand{\pnpo}{{\ensuremath{p^{n+1}}}}
\newcommand{\SGr}{\ensuremath{^{[S,G,r]}}}
\newcommand{\Sram}{\ensuremath{S_{\scriptstyle{ram}}}}

\newcommand{\td}{\ensuremath{\tilde{d}}}
\newcommand{\tDelta}{\ensuremath{{\tilde{\Delta}}}}

\newcommand{\tG}{\ensuremath{{\tilde{G}}}}
\newcommand{\tGamma}[1]{\ensuremath{\tilde{\Gamma}^{(#1)}}}
\newcommand{\tGinf}{\ensuremath{{\tilde{G}_{\infty}}}}
\newcommand{\tGinfd}{\ensuremath{{\tilde{G}^\dag_\infty}}}
\newcommand{\tH}{\ensuremath{\tilde{H}}}
\newcommand{\tHinf}{\ensuremath{{\tilde{H}_{\infty}}}}
\newcommand{\tiGamma}{\ensuremath{\tilde{\Gamma}}}
\newcommand{\tk}{\ensuremath{{\tilde{k}}}}
\newcommand{\tkinf}{\ensuremath{{\tilde{k}_{\infty}}}}

\newcommand{\tK}{\ensuremath{{\tilde{K}}}}
\newcommand{\tKinf}{\ensuremath{{\tilde{K}_{\infty}}}}
\newcommand{\tKn}{\ensuremath{{\tilde{K}_n}}}
\newcommand{\tL}[1]{\ensuremath{\tilde{L}^{(#1)}}}
\newcommand{\tpsi}{\ensuremath{\tilde{\psi}}}
\newcommand{\trho}{\ensuremath{\tilde{\rho}}}

\newcommand{\ttau}{\ensuremath{\tilde{\tau}}}
\newcommand{\talpha}{\ensuremath{\tilde{\alpha}}}

\newcommand{\tub}{\ensuremath{\underline{\tilde{b}}}}
\newcommand{\ub}{\ensuremath{\underline{b}}}
\newcommand{\ubi}{\ensuremath{\ub^{(i)}}}
\newcommand{\ueps}{\ensuremath{\underline{\eps}}}
\newcommand{\ueta}{\ensuremath{\underline{\eta}}}
\newcommand{\uu}{\ensuremath{\underline{u}}}
\newcommand{\uv}{\ensuremath{\underline{v}}}
\newcommand{\WedQGEinf}{\ensuremath{\bigwedge_{\bbQ \Ginf}^r\bbQ\cE_\infty(K)}}
\newcommand{\wmz}{\ensuremath{\fw_\fm^0}}
\newcommand{\ZGdd}{\ensuremath{\bbZ\times\Gammadd}}
\newcommand{\ZZnt}{\ensuremath{\bbZ\times \bbZ_{\geq n_2}}}
\begin{document}
\title{Abelian Conjectures of Stark Type\\
in $\bbZp$-Extensions of Totally Real Fields}
\author{D. Solomon
\thanks{partially supported by an Advanced Fellowship from the EPSRC.}
\\King's College London}
\maketitle
\section{Introduction}
We continue our investigations into complex and
$p$-adic variants of H.~M.~Stark's conjectures~\cite{Stark} for
an abelian extension of number fields $K/k$.
We have formulated versions of these conjectures
at $s=1$ using so-called `twisted zeta-functions'
(attached to additive characters) to replace the more usual $L$-functions.
The complex version of the conjecture was given
in~\cite{twizas}. In~\cite{zetap1} we formulated
an analogous $p$-adic conjecture in the case where $K$
is a (totally) real ray-class field over $k$,
as well as a third `combined' version that related the two previous ones
in this case.
(For numerical verification of the combined conjecture, see~\cite{zetap2}).
All these conjectures are stated in terms of group-ring-valued
regulators similar to those used in~\cite{Rubin}
to formulate fine `integral' variants of Stark's conjectures at $s=0$.

In Section~\ref{basic} of the present paper we restate
all three versions of the conjectures
mentioned above (designated respectively
$C1$, $C2$ and $C3$) in the context of an
essentially arbitrary, real abelian extension
$K/k$. However, we give
only the most `basic' forms of these conjectures \ie\
without any integrality conditions on their solutions, since
these conditions are of little direct concern to us in the present paper.
(See, however, Remarks~\ref{rem:denoms} and~\ref{rem:P1vsP1hat}).
Instead, our aim -- introduced in Section~\ref{section:Zpexts} --
is to focus on a new aspect
of Stark's Conjectures: under suitable conditions
both $C1$ and $C2$ predict the existence of `special'
elements in the $r$th exterior power of $\bbQ\otimes_{\bbZ}E(K)$ over
the group-ring $\bbQ \Gal(K/k)$ (where $r:=[k:\bbQ]$ and $E(K)$
denotes the unit group of $K$).
We consider the possible behaviour of these elements
as $K$ varies in a cyclotomic $\bbZ_p$-tower
$K_\infty=\bigcup_{n=0}^\infty K_n$.
More specifically, under the additional hypothesis that $p$ is unramified in $k$
we formulate a property `$P1$' which requires,
roughly speaking, that there exist solutions
of the complex conjecture $C1$ for
all $K_n/k$, $n\geq 1$, \emph{all coming from a single
element $\ueta$} of a module $\bigwedge^r\bbQ \cE_\infty$,
which is the exterior power over $\bbQ \Gal(K_\infty/k)$ of
the \emph{norm-coherent sequences} of (global) units
in the tower (tensored with $\bbQ$).
Property $P1$ is shown to hold in a
large infinite class of cases with $K\subset \bbQ^{\rm ab}$.
(Arguments given in Section~\ref{section:Zpexts}
make it seem plausible that $P1$ should hold rather
more generally. However, this is hard to verify, even numerically,
since doing so would seem to require
explicitly verifying $C1$ throughout the tower).
We also enunciate
an analogous property `$P2$' concerning the existence of such an $\ueta$
giving the solutions
of the $p$-adic conjecture $C2$ for the $\bbZ_p$-extension
(same $p$!) that is, for
each $K_n/k$ with $n$ larger than an explicit constant $n_2$.
We show that the
$p$-adic Property $P2$ follows from the complex one
$P1$ (with the same $\ueta$) provided we assume also the combined conjecture $C3$
for all $K_n/k$ with $n\geq n_2$.

Apart from their intrinsic interest and simplicity, the main
motivation for introducing properties $P1$ and, especially, $P2$
is revealed in Section~\ref{section:consP2}.
First, under the strengthened hypothesis that \emph{$p$ splits
in $k$}, we construct a doubly infinite
sequence of `higher' $p$-adic regulators $\check{\cR}_{t,n}$ on
$\bigwedge ^r\bbQ \cE_\infty$ by replacing
the  ($p$-adic) logarithm in the usual group-ring-valued
$p$-adic regulator with twisted,
generalised versions of the Coates-Wiles homomorphisms.
We then consider a putative
identity between two elements of the group-ring
$\bbC_p\Gal(K_n/k)$ for some given $n\geq n_2$ and $m\in\bbZ$.
The first element comes
from the $p$-adic twisted zeta-functions of the extension
$K_n/k$ evaluated at
$s=m$. The second element is
(up to an explicit algebraic factor)
the value of $\check{\cR}_{1-m,n}$
at a hypothetical element $\ueta$ of $\bigwedge^r\bbQ \cE_\infty$.
Thus it is, essentially, a determinant of `$(1-m)$th logarithmic
derivatives at level $n$' of norm-coherent sequences of global units.
The existence of an $\ueta$ satisfying
this relation for $m=1$ and all $n\geq n_2$
amounts to our $p$-adic Stark Conjecture $C2$
for each $K_n/k$ with $n\geq n_2$ as strengthened by saying precisely that
$\ueta$ satisfies
Property $P2$. On the other hand,
Theorem~\ref{thm:3A}
(the main result of this paper)
shows that if $\ueta$ satisfies our relation for infinitely many pairs
$(m,n)$ with $m\in\bbZ$ and $n\geq n_2$,
then it satisfies it for all such pairs.
The immediate consequence is, of course that if $\ueta$ satisfies
$P2$ then it also satisfies a new, `special value
conjecture' of Stark type for $p$-adic twisted
zeta-functions \emph{at any integral $s$}
(and this for each extension $K_n/k$ with $n\geq n_2$).

In the case $k=\bbQ$ of course, a 
link between the
Leopoldt-Kubota $p$-adic $L$-function and 
norm-coherent sequences of cyclotomic units was
already observed by Iwasawa in~\cite{Iwasawa}. Moreover the use of
the Coates-Wiles homomorphisms to make this connection is now well-known
(see \eg~\cite{Wa}). One way to view the present paper
is as an explanation of how
a strengthened $p$-adic Stark Conjecture at $s=1$ leads to a 
related and rather precise connection for twisted zeta-functions 
and more general $K/k$ as above. 
On the other hand, a theory due to 
Perrin-Riou proposes extremely general 
links between norm-coherent sequences 
in the first cohomology of a $p$-adic representation of 
$\Gal(\bar{\bbQ}/\bbQ)$ and an associated 
$p$-adic $L$-function (whose existence is partly 
conjectural. See~\cite{PR2} for details).
If this latter theory, or an extension of it, 
could be made to produce explicit predictions
for $L$-functions of the extension $K/k$, then interesting
comparisons with our conjectures might follow.
So far this seems only to have been done in the above known case
of the conjectures, 
namely $k=\bbQ$ (see \eg~\cite{PR1} and references therein).

In section~\ref{section:pfmainthm} we give the
proof of Theorem~\ref{thm:3A} which relies at base on
a well known uniqueness
principle in $p$-adic analysis stating roughly that a non-zero
$p$-adic power series with bounded
coefficients can have only finitely many roots in $\bbC_p$.
To apply this principle we use
the theory of $p$-adic measures and results of Deligne and
Ribet on the analytic behaviour of $p$-adic $L$-functions attached to
characters of $\Gal(K_n/k)$. Global Gauss sums
intervene in relating the $p$-adic $L$-functions to our $p$-adic twisted
zeta-functions, so the proof also requires a rather intricate
comparison between these sums and a certain product of
local Gauss sums that appears naturally elsewhere.
(It should be possible to avoid all Gauss sums by using the twisted
zeta-functions directly. However, since their $p$-adic behaviour is less
well-understood that of $L$-functions, we have so far done
this only for certain classes of real-quadratic $k$, using
different methods coming from~\cite{plim}.)

Finally, in Section~\ref{section:semiloc} we show how
$P1$ and $P2$ can be weakened by `enlarging'
$\bigwedge^r\bbQ \cE_\infty$ to
$\bbQ_p\otimes_{\bbZ_p}\bigwedge^r\cU^{(1)}_\infty$.
where, $\cU^{(1)}_\infty$ is the module of
norm-coherent sequences of
\emph{$p$-semilocal} units in the tower and the exterior power is taken over
the appropriate \emph{completed $p$-adic group-ring}.
Thus, for example, $\widehat{P2}$ weakens
$P2$ by assuming only that
there exists an element
$\hat{\ueta}\in\bbQ_p\otimes_{\bbZ_p}\bigwedge^r\cU^{(1)}_\infty$
whose projection to the $n$th level is in the module
of exterior powers of \emph{global} units of $K_n$ for all $n$
and that this projection
is also a solution of $C2$ for $K_n/k$ if $n\geq n_2$.
The higher regulators $\check{\cR}_{t,n}$
`extend' to $\bbQ_p\otimes_{\bbZ_p}\bigwedge^r\cU^{(1)}_\infty$ and the same
methods as in the preceding sections
show that the weaker existence hypothesis
$\widehat{P2}$ implies the weaker conclusion that there is exists an element
$\hat{\ueta}\in\bbQ_p\otimes_{\bbZ_p}\bigwedge^r\cU^{(1)}_\infty$
which is global at each finite level and
satisfies the special value conjectures at all $s\in\bbZ$ and $n\geq n_2$
(namely, $\hat{\ueta}$ is the element satisfying $\widehat{P2}$).

{\bf Some Notation:} For the rest of this paper we fix the totally real number field $k$
considered as contained in $\barbbQ\subset\bbC$. In addition to
the above notations, we shall write
$\cO$ for the ring of integers of $k$, $d_k\in\bbN$ for its absolute discriminant
and $S_\infty$ for the set of its infinite places.
The latter may be identified with the (real) embeddings
of $k$ into $\barbbQ$ which we denote $\tau_1,\ldots,\tau_r$
and for each $i=1,\ldots,r$ we fix once and for all
an element $\tilde{\tau}_i\in\GQbarQ$ extending $\tau_i$.
For any $m\in\bbZ_{\geq 1}$ we write
$\mu_m$ for the group of all $m$th roots of unity in $\barbbQ$
and $\zeta_m$ for its specific generator $\exp(2\pi i/m)$.

\section{Basic Conjectures}\label{basic}
We start by restating Conjectures~3.1, 3.2 and 3.3
of~\cite{zetap1}, referring to~\emph{ibid.} for many of the
details. First, recall the definitions
of the functions $\PhimT$ and $\PhimTp$.
Let $\fm$ be a cycle for $k$, namely a formal product $\fg\fz$
where $\fg$ is a non-zero ideal of $\cO$ and $\fz$ is (the formal
product of) a subset of $S_\infty$. The ray-class
group and the ray-class field of $k$ modulo $\fm$ will be denoted
$\Cl_\fm(k)$ and $k(\fm)$ respectively (thus $k\subset k(\fm)\subset\barbbQ$).
Let $T$ be any finite subset
of prime ideals of $\cO$. We defined in \cite{twizas,zetap1} a
function of one complex variable $s$ taking values in the
group-ring $\bbC \Gal(k(\fm)/k)$ by setting
\[
\PhimT(s)=\sum_{\fc\in\Clmk}Z_T(s;\fc\cdot\wmz)\sigma_\fc\inv
\]
Here $\sigma_\fc$ is the element of $\Gal(k(\fm)/k)$ associated
to $\fc$ by the Artin isomorphism and `$Z_T(s;\fc\cdot\wmz)$' is a
\emph{twisted zeta-function} depending on $\fc$.
We emphasize that in the present context we are \emph{not}
assuming --- as has been done in some similar ones --- that
$\fg$ is `prime to $T$', \ie\ that $v_\fp(\fg)=0\ \forall \fp\in T$.
We shall say that $\fg$ is \emph{$T$-trivial}
if and only if it is a product --- possibly empty ---
of \emph{distinct} primes of $T$.
Recall that $Z_T(s;\fc\cdot\wmz)$ was defined by a
Dirichlet series for ${\rm Re}(s)>1$
and that it (and hence also $\PhimT(s)$)
is closely related to the $L$-functions
of the characters of $\Clmk$
(see~Proposition~\ref{prop:1A}).
Recall (\cite[Thm. 2.3]{zetap1})
that $\PhimT(s)$ extends to a meromorphic function on $\bbC$
with at most a simple pole
at $s=1$ and is holomorphic in $\bbC$ iff $\fg$ is not $T$-trivial.

The values of $\PhimT$ at non-positive integers lie in the group
ring $\barbbQ \Gal(k(\fm)/k)$ (see~\cite[Lemma 3.2]{zetap1}). Under
certain conditions they can be $p$-adically interpolated. Let $p$
be any prime number. We shall use the following notation.
Every $x\in\bbZ_p^\times$ can be written uniquely as
\[
x=\omega(x)\langle x\rangle
\]
where $\omega(x)$ is a root of unity in $\bbZ_p^\times$ and
$\langle x\rangle$ lies in $1+p\bbZ_p$ (in $1+4\bbZ_2$ if $p=2$).
Let
$\cM(p)$ denote the set of integers $m\in\bbZ_{\leq 0}$ such that
$m\equiv 1 \pmod{p-1}$ ($({\rm mod}\ 2)$ if $p=2$).
Let $D(p)$ denote the $p$-adic disc $1+2\bbZp$ and
$D^0(p)$ the punctured disc $D(p)\setminus\{1\}$.
Let $\bbCp$ be a completion of an algebraic
closure $\barbbQ_p$ of $\bbQ_p$. Now fix an embedding
$j:\barbbQ\rightarrow\bbCp$ and let it act on coefficients to
define a homomorphism from $\barbbQ \Gal(k(\fm)/k)$ to $\bbC_p
\Gal(k(\fm)/k)$. We showed
(see~\cite[Thm./Def. 2.1]{zetap1}) that if $T$ contains the set
$T_p$ of places of $k$ dividing $p$, then then there exists a
unique, $p$-adically continuous function $\PhimTp=\PhijmTp$ from
$D^0(p)$ to $\bbCp \Gal(k(\fm)/k)$ such that
$\PhimTp(m)=j(\PhimT(m))$ for every $m\in\cM(p)$. In fact, $\PhimTp$ is
meromorphic on $D(p)$ and even holomorphic there if
$\fg$ is not $T$-trivial. It is related to
the \emph{$p$-adic} $L$-functions which can be attached to the
(totally) even characters of $\Clmk$ and consequently
depends only on $\fg$, not on the infinite part $\fz$
of $\fm$. (See Proposition~\ref{prop:1A} for a more precise statement).

The `basic' conjectures of~\cite{zetap1} concern the values $\PhimT(1)$ and
$\PhimTp(1)$ in the case where $\fz$ is trivial so that
$k(\fm)=k(\fg)$ is a finite \emph{real} abelian extension of $k$.
Generally, we shall write $K$ for such an extension (as in the introduction)
and $G$ for $\Gal(K/k)$.
Given also a finite set $S$
of places $k$ containing $S_\infty$, we
write $U_S(K)$ for the group of $S$-units of $K$.
We define $\bbZ G$-linear,
logarithmic maps
$\lambda_i=\lambda_{K/k,S,i}:U_S(K)\rightarrow \bbR G$ and
$\lambda_{i,p}=\lambda_{K/k,S,i,p}^{(j)}:U_S(K)\rightarrow \bbCp G$ by
setting $\lambda_i(u)= {\sum_{\sigma\in
G}}\log|\tilde{\tau}_i\sigma(u)|\sigma\inv$ and
$\lambda_{i,p}(u)={\sum_{\sigma\in G}}
\log_p(j\tilde{\tau}_i\sigma(u))\sigma\inv$ for any $u\in U_S(K)$
where $\log_p$ denotes Iwasawa's $p$-adic logarithm.
Each such map `extends' by $\bbQ$-linearity to
$\bbQ U_S(K):=\bbQ\otimes_\bbZ U_S(K)$.
Thus we can define unique, $\bbQ G$-linear `regulator' maps
$R_{K/k}=R_{K/k,S}$
and $R_{K/k,p}=R_{K/k,S,p}^{(j)}$ from the $r$th exterior power
$\bigwedge_{\bbQ G}^r\bbQ U_S(K)$ (notated as an additive
$\bbQ G$-module) into $\bbR G$
and $\bbC_p G$ respectively by setting $R_{K/k}
(u_1\wedge\ldots\wedge u_r)=\det(\lambda_i(u_t))_{i,t=1}^r$ and
$R_{K/k,p} (u_1\wedge\ldots\wedge u_r)=
\det(\lambda_{i,p}(u_t))_{i,t=1}^r$ for all
$u_1,\ldots,u_r\in\bbQ U_S(K)$.

In the case $K=k(\fg)$, we write
$G_\fg$ for $\Gal(k(\fg)/k)$ and let $S(\fg)$ be the set of places of $k$
dividing $\fg$ together with the infinite ones.
Conjecture~3.1
of~\cite{zetap1} then reads as follows
\begin{conjalt}{$C1(k,\fg,T)$} Suppose that $\fg$ is not $T$-trivial.
Then
\beql{eq:1A} \Phi_{\fg,T}(1)=
\frac{2^r}{\sqrt{d_k}\prod_{\fp\in T}N\fp} R_{k(\fg)/k}(\eta)\ \ \
\mbox{for some $\eta\in\bigwedge_{\bbQ G_\fg}^r\bbQ U_{S(\fg)}(k(\fg))$}
\eeq
\end{conjalt}
\noindent Note that the conditions of
this conjecture already imply that $\Phi_{\fg,T}(1)\in\bbR G_\fg$
(see Lemma~3.1 of~\cite{twizas}). A $p$-adic analogue
(Conjecture~3.2 of~\cite{zetap1})
can be formulated as
\begin{conjalt}{$C2(k,\fg,T,p)$}
Suppose that $T$ contains $T_p$ and that $\fg$ is not $T$-trivial. Then
\beql{eq:1B} \Phi_{\fg,T,p}(1)=
\frac{2^r}{j(\sqrt{d_k})\prod_{\fp\in T}N\fp}
R_{k(\fg)/k,p}(\eta)\ \ \
\mbox{for some $\eta\in\bigwedge_{\bbQ G_\fg}^r\bbQ U_{S(\fg)}(k(\fg))$} \eeq
\end{conjalt}
\noindent Finally, we also made a
combined conjecture (Conjecture~3.3 in~\cite{zetap1}):
\begin{conjalt}{$C3(k,\fg,T,p)$} Under the conditions of $C2(k,\fg,T,p)$ there
exists a simultaneous solution
$\eta\in\bigwedge_{\bbQ G_\fg}^r\bbQ U_{S(\fg)}(k(\fg))$
of both\refeq{1A} and\refeq{1B}.
\end{conjalt}
\noindent In\refeq{1A} and\refeq{1B}
the rational normalising factor $2^r/(\prod_{\fp\in
T}N\fp)$ could of course be absorbed into $\eta$
but it is more convenient to keep it. The dependence of $\eta$ on $T$ and the choices
of the $\tilde{\tau_i}$ --- as well as its indepedence of $j$
in\refeq{1B} --- are explained in \S~3.2 of~\cite{zetap1}.
In \S~3.3 \emph{ibid.} we examined the relations
between these conjectures and
Stark's original conjectures as well as the complex and
$p$-adic variants
given by the author, Rubin, Serre \emph{etc.}
(see \cite{twizas},~\cite{Rubin},~\cite{Tate}\ldots).
A refined version of Conj. $C3(k,\fg,T,p)$
was formulated as Conjecture~3.6 of~\cite{zetap1}.
This was the conjecture tested numerically in~\cite{zetap2}
(always with $k$ quadratic and $\fg$ prime to
$T=T_p$).

Now let $K$ be a general finite, real, abelian extension of $k$
with group $G$.
The above conjectures generalise easily from the case  $K=k(\fg)$
by `taking norms from the conductor field'.
More precisely, if $K'$
is any other such extension
containing $K$, with $G':=\Gal(K'/k)$, we shall write
$\pi_{K'/K}$ for the restriction $G'\rightarrow G$,
linearly extended to a homomorphism of group-rings.
The norm homomorphism $N_{K'/K}$ defines a $\pi_{K'/K}$-semilinear map
from $U_S(K')$ into $U_S(K)$ for any $S$, so it induces
a homomorphism (also denoted
$N_{K'/K}$) from $\bigwedge_{\bbQ G'}^r\bbQ U_S(K')$ to
$\bigwedge_{\bbQ G}^r\bbQ U_S(K)$,
for which one checks easily that
\beql{eq:1B.5}
R_{K/k}\circ N_{K'/K}=\pi_{K'/K}\circ R_{K'/k}
\eeq
The conductor $\ff(K)$ of $K$ is the largest
ideal $\ff$ of $\cO$ such
that $k\subset K\subset k(\ff)$.
It can be determined by local analysis of ramification
(see\refeq{2F}).
In particular, $S(\ff(K))$ equals $S_{\scriptstyle{ram}}=
S_{\scriptstyle{ram}}(K/k)$
which we define to be the set of finite places ramified
in $K/k$ (which are also those ramified in $k(\ff(K))/K$)
\emph{together with $S_\infty$}.
Let $\cD=\cD_{k/\bbQ}$ denote the absolute different.
For any $T$ we denote by $\ff'(K)$ and $\cD'$ the prime-to-$T$ parts of
$\ff(K)$ and $\cD$ respectively (in the obvious sense)
and
define a function $\PhiKkT:\bbC\setminus\{1\}\rightarrow \bbC G$ by setting
\beql{eq:1C}
\PhiKkT(s):=N(\ff'(K)\cD')^{s-1}\pi_{k(\ff(K))/K}(\Phi_{\ff(K),T}(s))
\eeq
If $p$ is a prime number and $T$ contains $T_p$ then
there exists a $p$-adically continuous function
$\PhiKkTp=\Phi^{(j)}_{K/k,T,p}:D^0(p)\rightarrow \bbCp G$
such that
$\PhiKkTp(m)=j(\PhiKkT(m))$ for every $m\in\cM(p)$.
Indeed such a function
is clearly unique and given by
\beql{eq:1D}
\PhiKkTp(s):=\langle N(\ff'(K)\cD')\rangle ^{s-1}
\pi_{k(\ff(K))/K}(\Phi_{\ff(K),T,p}(s))
\eeq
If $\ff(K)$ is not $T$-trivial then $\PhiKkT$ and
$\PhiKkTp$ extend to holomorphic
functions on $\bbC$ and $D(p)$ respectively.\vspace{1ex}\\
\rem\
In~\cite{twizas} a
complex function $\Phi^{-}_{K,\fm,T}$ rather more general
than $\PhiKkT$ was defined by a similar process,
but without the factor $N(\ff'(K)\cD')^{s-1}$. Of course, neither
this nor the factor
$\langle N(\ff'(K)\cD')\rangle ^{s-1}$ in\refeq{1D}
has any effect on conjectures about values at $s=1$. Both factors are
introduced here for purposes of `normalisation' and for conjectures at other
integral values of $s$.
\vspace{1ex}\\
Now take $K'=k(\ff(K))$ and apply $\pi_{k(\ff(K))/K}$
to\refeq{1A} and\refeq{1B} with $\fg=\ff(K)$.
Using\refeq{1B.5}
(with $K'=k(\ff(K))$) and replacing $N_{k(\ff(K))/K}\eta$ by $\eta$,
we see that Conjectures~$C1(k,\ff(K),T)$, $C2(k,\ff(K),T,p)$
and $C3(k,\ff(K),T,p)$ imply respectively
\begin{conjalt}{$C1(K/k,T)$} Suppose that $\ff(K)$ is not $T$-trivial.
Then
\beql{eq:1E} \PhiKkT(1)=
\frac{2^r}{\sqrt{d_k}\prod_{\fp\in T}N\fp} R_{K/k}(\eta)\ \ \
\mbox{for some $\eta\in\bigwedge_{\bbQ G}^r\bbQ U_{\Sram}(K)$}
\eeq
\end{conjalt}
\begin{conjalt}{$C2(K/k,T,p)$} Suppose that $T$ contains $T_p$ and
$\ff(K)$ is not $T$-trivial. Then
\beql{eq:1F} \PhiKkTp(1)=
\frac{2^r}{j(\sqrt{d_k})\prod_{\fp\in T}N\fp} R_{K/k,p}(\eta)\ \ \
\mbox{for some $\eta\in\bigwedge_{\bbQ G}^r\bbQ U_{\Sram}(K)$}
\eeq
\end{conjalt}
\begin{conjalt}{$C3(K/k,T,p)$} Under the conditions of $C2(K/k,T,p)$ there
exists a simultaneous solution
$\eta\in\bigwedge_{\bbQ G}^r\bbQ U_{\Sram}(K)$
of both\refeq{1E} and\refeq{1F}.
\end{conjalt}
\noindent It will sometimes be convenient to
use the phrase `a solution of $C1(K/k,T)$'
(resp. of $C2(K/k,T,p)$', resp. of $C3(K/k,T,p)$')
to mean an element $\eta$ of
$\bigwedge_{\bbQ G}^r\bbQ U_{\Sram}(K)$ satisfying\refeq{1E}
(resp.\refeq{1F}, resp.\refeq{1E} and\refeq{1F}).\vspace{1ex}\\
\rem\
If $\fg$ is any non-$T$-trivial ideal,
then $C1(k(\fg)/k,T)$ implies $C1(k,\fg,T)$
\emph{provided} that the ideal $\ff(k(\fg))$ is also non-$T$-trivial.
(This is the case, for instance, if $\ff(k(\fg))=\fg$, \ie\ if
$\fg$ is a conductor).
Indeed if we set $K=k(\fg)$ then $\ff(K)$ divides
$\fg$ and $K$ also equals $k(\ff(K))$.
Theorem~2.1 of~\cite{zetap1} then implies that $\PhigT(1)$
is a multiple of $\Phi_{\ff(K),T}(1)=\Phi_{K/k,T}(1)$
in $\bbC Gal(K/k)$ so that\refeq{1E} implies\refeq{1A}.
Similar arguments work for $C2$ and $C3$ with similar restrictions.
\vspace{1ex}\\
We record some basic notations and results to be used in
investigating these conjectures.
\begin{lemma}\label{lemma:1Z}
Suppose that $K'\supset K$ as above and (for simplicity) that $\Sram(K'/k)$ equals
$\Sram(K/k)$ and is disjoint from $T$. Then
\[
\pi_{K'/K}(\Phi_{K'/k,T}(s))=
\PhiKkT(s)
\]
\end{lemma}
\bPf\
$\ff(K')$ has the same prime factors as
$\ff(K)$ and is prime to $T$, so
Theorem~3.2 and Remark~3.1 of~\cite{twizas} give
\[
\pi_{k(\ff(K'))/k(\ff(K))}(\Phi_{\ff(K'),T}(s))=
\left(\frac{N\ff(K')}{N\ff(K)}\right)^{1-s}
\Phi_{\ff(K),T}(s)
\]
Now apply $\pi_{k(\ff(K))/K}$ to
both sides and use $\pi_{k(\ff(K))/K}\circ\pi_{k(\ff(K'))/k(\ff(K))}=
\pi_{K'/K}\circ\pi_{k(\ff(K'))/K'}$ and\refeq{1C} (twice).\ePf
\noindent For any profinite abelian group $A$ we write $A^\dag$ for
the group of (degree-$1$) complex characters of $A$ which factor through some finite
quotient of $A$.
Any $\chi\in G^\dag$ gives rise to a $p$-adic character by
composition with $j$ which, by abuse, will also be denoted
$\chi$. So also will the associated characters (complex or
$p$-adic) of $\Gal(k(\ff(K))/k)$ and $\Cl_{\ff(K)}(k)$, obtained from $\chi$
\emph{via} $\pi_{k(\ff(K))/K}$ and the Artin map respectively.
The conductor $\ff(\chi)$ of $\chi$
is that of the subextension $K^{\ker(\chi)}/k$ that it cuts out. Thus
$\ff(K)=\fh\ff(\chi)$ for some integral ideal $\fh$. We write
$\ff'(\chi)$ and $\fh'$ for the prime-to-$T$ parts so that
$\ff'(K)=\fh'\ff'(\chi)$ and and $\fh=\fh'\fh_0$ where $\fh_0$
has support in $T$. We write $\hchi$ for the unique, complex or
$p$-adic primitive character on $\Cl_{\ff(\chi)}(k)$ associated
to $\chi$ and $\hchi(\fa)$ for its value on the class of an ideal
$\fa$ which is prime to $\ff(\chi)$. If $\chi$ is complex then
the Gauss sum $g_{\ff(\chi)}(\hchi)\in\barbbQ^\times$
defined in~\cite[\S 6.4]{twizas}
will  be denoted simply $g(\hchi)$. If $\chi$ is $p$-adic, $\chi=j\circ\phi$
then we define $g(\hchi)\in\bbC_p$ to be $j(g(\hat{\phi}))$.
(Note that this is \emph{not} independent of the choice of $j$.)
With these conventions the following result relates
$\PhiKkT(s)$ and $\PhiKkTp(s)$ to the (primitive) complex and
$p$-adic $L$-functions, denoted $L(s,\hchi)$ and $L_p(s,\hchi)$
respectively.
\begin{prop}\label{prop:1A}
Suppose that $\chi$ is complex and $s\in\bbC\setminus\{1\}$,
\emph{respectively} $\chi$ is $p$-adic, $T\supset T_p$ and
$s\in D^0(p)$. If $\fh_0$ is a product of $t\geq 0$
distinct primes not dividing $\ff(\chi)$ then
\begin{eqnarray}
\lefteqn{\chi(\PhiKkT(s))=(-1)^t\hchi\inv(\fh_0)g(\hchi)(N(\ff'(\chi)\cD'))^{s-1}
\prod_{\fp\ \mathit{prime}\atop \fp|\fh',\ \fp\ndiv\ff(\chi)}
\left(1-\frac{\hchi\inv(\fp)}{N\fp^{1-s}}\right)\times}\hspace{24em}\nonumber\\
 \hspace{20em}&&\prod_{\fp\
\mathit{prime}\atop \fp\in T,\ \fp\ndiv \ff(\chi)}
\left(1-\frac{\hchi(\fp)}{N\fp^s}\right)L(s,\hchi)\nonumber\\
&&\label{eq:1G}
\end{eqnarray}
respectively
\begin{eqnarray}
\lefteqn{\chi(\PhiKkTp(s))=(-1)^t\hchi\inv(\fh_0)g(\hchi)
\langle N(\ff'(\chi)\cD')\rangle ^{s-1}
\prod_{\fp\ \mathit{prime}\atop \fp|\fh',\ \fp\ndiv\ff(\chi)}
\left(1-\frac{\hchi\inv(\fp)}{\langle N\fp\rangle ^{1-s}}\right)\times}\hspace{20em}
\nonumber\\
\hspace{20em}&&
\prod_{\fp\ \mathit{prime}\atop \fp\in T,\ \fp\ndiv p\ff(\chi)}
\left(1-\frac{\hchi(\fp)}{\omega(N\fp)\langle N\fp\rangle
^s}\right)
L_p(s,\hchi)\nonumber\\
&&\label{eq:1H}
\end{eqnarray}
Otherwise
\beql{eq:1I}
\chi(\PhiKkT(s))=0,
\eeq
respectively
\beql{eq:1J}
\chi(\PhiKkTp(s))=0
\eeq
\end{prop}
\bPf\ This follows easily by Equations\refeq{1C} and \refeq{1D}
from Theorem~2.2 and part~(iv) of Theorem/Definition~2.1 of~\cite{zetap1}.
\ePf
\noindent The behaviour of $L$-functions in a neighbourhood of
$s=1$ shows that the R.H.S. of\refeq{1G} and\refeq{1H} define
analytic functions there whenever $\chi$ is not
equal to the trivial character $\chi_0$
or $\ff(K)$ is not $T$-trivial.
In particular, taking limits as $s\rightarrow 1$, we obtain
\begin{cor}\label{cor:1A} If $\ff(K)$ is not $T$-trivial then
Equations\refeq{1G},\refeq{1H},\refeq{1I}
and\refeq{1J} also hold for $s=1$ and any character $\chi$ of $G$.
(For $\chi=\chi_0$,
we must take $\lim_{s\rightarrow 1}$ on the R.H.S. in\refeq{1G} and\refeq{1H}.)\ePf
\end{cor}
For any $\chi\in G^\dag$ we write $e_\chi$ for the idempotent
$|G|\inv\sum_{\sigma\in G}\chi(\sigma)\sigma\inv\in\bbC G$
and for any $S\supset S_\infty$ as above we set
\beql{eq:1K}
r(S,\chi):={\rm dim}_\bbC(e_\chi \bbC U_S(K))=\left\{
\begin{array}{ll}
r+|\{\mbox{prime ideals}\ \fq\in S: \chi|_{G(\fq)}=1\}|&
\mbox{if $\chi\neq\chi_0$}\\
r-1+|\{\mbox{prime ideals}\ \fq\in S\}|=|S|-1 & \mbox{if $\chi=\chi_0$}
\end{array}
\right.
\eeq
(Here $G(\fq)$ denotes the decomposition subgroup of
$G$ associated to $\fq$.
For the second equality, see~\cite[\S\S I.3, I.4]{Tate}.)
Then Proposition~\ref{prop:1A} shows that
\beql{eq:1L}
\mbox{$T$ is prime to $\ff(K)$}\Longrightarrow
\ord_{s=1}(\chi(\PhiKkT(s)))=r(\Sram(K),\chi)-r\ \forall\chi\in G^\dag
\eeq
If $S\neq S_\infty$  then $r(S,\chi)\geq
r\ \forall \chi\in G^\dag$.
In this case, we denote by $e_{S,G,r}$ (respectively $e_{S,G,>r}$) the
sum of the idempotents $e_\chi\in \bbC G$ for those $\chi$ with
$r(S,\chi)=r$ (respectively $r(S,\chi)>r$).
Clearly, $r(S,\chi)$ depends only on the
$\Gal(\barbbQ/\bbQ)$-orbit $[\chi]$ of $\chi$
(so we write also $r(S,[\chi])$).
Hence $e_{S,G,r}$ is a sum of $e_{[\chi]}$'s where
$e_{[\chi]}:=\sum_{\chi'\in[\chi]}e_{\chi'}\in\bbQ G$
and similarly for $e_{S,G,>r}$. In particular, $e_{S,G,r}$ and
$e_{S,G,>r}$ are complementary idempotents of $\bbQ G$ and
$\tilde{e}_{S,G,>r}:=|G|e_{S,G,>r}$ lies in $\bbZ G$. For any
$G$-module $A$ we let $A^{[S,G,r]}:=\ker \tilde{e}_{S,G,>r}|A$.
If $A$ is also a $\bbQ$-vector space (and $S\neq S_\infty$)
then $A^{[S,G,r]}=e_{S,G,r}A=\bigoplus_{r(S,[\chi])=r}e_{[\chi]}A$ on which
the $\bbQ G$-action factors through the quotient ring
$e_{S,G,r} \bbQ G=\prod_{r(S,[\chi])=r}F_{[\chi]}$. Here
$F_{[\chi]}$ denotes $e_{[\chi]}\bbQ G$ which is a \emph{field}
(isomorphic to $\bbQ(\chi)$ via $\chi$).
\begin{lemma}\label{lemma:1A} If $S\neq S_\infty$ then
$R_{K/k,S}$ is injective on $(\bigwedge_{\bbQ G}^r\bbQ U_S(K))^{[S,G,r]}$.
\end{lemma}
\bPf One can deduce this lemma from the existence of a solution of $C1(K/k,\emptyset)$
(\cf\ the proof of Proposition~3.8~(i) in~\cite{zetap1})
and this would suffice for our applications. However, it also follows unconditionally
from Remark~2, \S 1.6 of~\cite{Po1}. Indeed, the module
$\bbC\otimes_\bbQ(\bigwedge_{\bbQ G}^r\bbQ U_S(K))^{[S,G,r]}$ in our notation
coincides with Popescu's `$(\bbC\bigwedge^r U_{S,T})_{r,S}$'
(the set $T$ is irrelevant here and may be chosen conveniently).
Moreover, the $\bbC$-linear extension of our map $R_{K/k,S}$
coincides (up to sign) with Popescu's `$R_W$' on $\bbC\bigwedge^r U_{S,T}$
provided we take his `$w_i$' to be the places defined by our
$\ttau_i$, for $i=1,\ldots,r$. Popescu's context is function-fields
but his proof goes over to number fields without real change.
Note also that the Lemma actually holds without the assumption $S\neq S_\infty$.
\ePf\noindent
As we shall see later, if a solution of\refeq{1E} exists
then certain conditions imply that one exists in
the subspace
$(\bigwedge_{\bbQ G}^r\bbQ U_{\Sram}(K))^{[\Sram,G,r]}$ and,
using the above lemma,
that such a solution is unique.
The other interest of this subspace comes from the fact that every element is
a $\bbQ$-multiple of a wedge product of actual units of $K$
(provided that
$\Sram$ contains at least two finite places). More precisely, if
$\varepsilon\in E(K)=U_{S_\infty}(K)$
is such a unit then we shall write $\bareps$ for its
image $1\otimes \eps$ in $\bbQ E(K)\subset \bbQ U_S(K)$.
The kernel of the map $\eps\mapsto\bareps$ is $E(K)_{\rm tor}=\{\pm 1\}$
and we have
\begin{lemma}\label{lemma:1B}
Suppose that $S$ contains at least two finite places. Then
every element of $(\bigwedge_{\bbQ G}^r\bbQ U_S(K))^{[S,G,r]}$
can be written
$\frac{1}{c}\bar{\varepsilon}_1\wedge\ldots\wedge\bar{\varepsilon}_r$ for
some $c\in\bbN$ and $\varepsilon_1,\ldots,\varepsilon_r\in
E(K)^{[S,G,r]}$.
\end{lemma}
\bPf\
Let $\chi$ be an element of $G^\dag$. It is easy
to show that $\dim_{F_{[\chi]}}(e_{[\chi]}\bbQ U_S(K))=
{\rm dim}_\bbC(e_\chi \bbC U_S(K))=r(S,\chi)$
and similarly that
$\dim_{F_{[\chi]}}(e_{[\chi]}\bbQ E(K))=r(S_\infty,\chi)$.
But the condition on $S$ shows that if $r(S,\chi)=r$
then $\chi\neq\chi_0$ so $r(S_\infty,\chi)$ also equals $r$.
Hence $e_{[\chi]}\bbQ U_S(K)=e_{[\chi]}\bbQ E(K)$
for all such $[\chi]$. Summing, we deduce that $(\bbQ U_S(K))^{[S,G,r]}$
equals $(\bbQ E(K))^{[S,G,r]}=\bbQ (E(K)^{[S,G,r]})$
and is free of rank $r$ over
$\prod_{[\chi]}F_{[\chi]}=\bbQ G\SGr$.
Consequently, every element of
$(\bigwedge_{\bbQ G}^r\bbQ U_S(K))\SGr
=\bigwedge_{\bbQ G\SGr}^r((\bbQ U_S(K))^{[S,G,r]})$ can be written as
$u_1\wedge\ldots\wedge u_r$ with
$u_1,\ldots,u_r\in\bbQ (E(K)^{[S,G,r]})$.\ePf

\section{$\bbZ_p$-Extensions}\label{section:Zpexts}
For the rest of this paper we fix an extension $K/k$ with
Galois group $G$, a prime number $p$ and an embedding
$j:\barbbQ\rightarrow\bbC_p$ all subject to the
hypotheses of the previous section
as well as the following conditions
\beql{eq:2A} p\neq 2
\eeq
\beql{eq:2B} \mbox{$\Sram(K/k)$ contains at
least one finite place not dividing $p$}
\eeq
and
\beql{eq:2C} \mbox{$p$ is unramified in $k/\bbQ$}
\eeq
Condition\refeq{2A} will be convenient but
not necessary in all that follows. Conditions\refeq{2B}
and \refeq{2C} could possibly
be weakened for this section but\refeq{2C}
anyway will need to be strengthened for the next.

For $n=0,1,2,\ldots$ let $\bbB_n$ be the unique subfield
of $\bbQ(\mu_{p^{n+1}})$ which is cyclic and of degree $p^n$ over $\bbQ$.
For any abelian extension $L$ of $k$ we write $L_n$ for $L\bbB_n$
(so $L_0=L$). Thus $K_n$ is a totally real
abelian extension of $k$ with group $G_n$ say.
It is cyclic over $K$ of degree \emph{dividing} $p^n$.
We let $K_\infty$ be $\bigcup_{n\geq 0}K_n$, the cyclotomic $\bbZ_p$-extension of $K$
and write $G_\infty$ for $\Gal(K_\infty/k)$.
Let $n_1=n_1(K,p)$ be the smallest integer $n$ such that
$K_n/\bbQ$ is ramified at all primes above $p$. Thus $n_1=0$ or $1$ and
Condition\refeq{2C} implies that the set
$\Sram(K_n/k)$ equals $\Sigma=\Sigma(K/k,p):=\Sram(K/k)\cup T_p$
for all $n\geq n_1$ hence it
contains at least two finite places by Condition\refeq{2B}.
In particular, $\ff(K_n)$ is not $\emptyset$-trivial,
so $\Phi_{K_n/k,\emptyset}$ is holomorphic and $C1(K_n/k,\emptyset)$
makes sense.
\begin{lemma}\label{lemma:2A}
Suppose $n\geq n_1$ and $\eta$ is a solution of
$C1(K_n/k,\emptyset)$ then $e_{\Sigma,G_n,r}\eta$
is the unique solution lying in
$(\bigwedge_{\bbQ G_n}^r\bbQ U_{\Sigma}(K_n))^{[\Sigma,G_n,r]}$
\end{lemma}
\bPf\
The condition that $\eta$ be a solution of
$C1(K_n/k,\emptyset)$ is clearly equivalent to
\beql{eq:2C.5}
\frac{2^r}{\sqrt{d_k}}\chi(R_{K_n/k}(\eta))=\chi(\Phi_{K_n/k,\emptyset}(1))
\ \ \ \forall \chi\in G_n^\dag
\eeq
To show that $e_{\Sigma,G_n,r}\eta$ is also a solution,
we must deduce that\refeq{2C.5}
also holds with $e_{\Sigma,G_n,r}\eta$ in place of $\eta$. But
$\chi(R_{K_n/k}(e_{\Sigma,G_n,r}\eta))$ is equal to
$\chi(e_{\Sigma,G_n,r})\chi(R_{K_n/k}(\eta))$
and hence to $\chi(R_{K_n/k}(\eta))$
or $0$ according as $r(\Sigma,\chi)=r$ or $r(\Sigma,\chi)>r$,
and in the latter case,\refeq{1L} shows that one also has
$\chi(\Phi_{K_n/k,\emptyset}(1))=0$, as required.
It is clear that $e_{\Sigma,G_n,r}\eta$ lies
in $(\bigwedge_{\bbQ G_n}^r\bbQ U_{\Sigma}(K_n))^{[\Sigma,G_n,r]}$
and the uniqueness follows from Lemma~\ref{lemma:1A}.\ePf
\noindent
We shall be concerned primarily with the
complex functions $\Phi_{K_n/k,\emptyset}$ and $\Phi_{K_n/k,T_p}$
(for $n\geq n_1$) as well as the $p$-adic function
$\Phi^{(j)}_{K_n/k,T_p,p}$ which interpolates $\Phi_{K_n/k,T_p}$.
When no confusion is possible we shall abbreviate these three functions respectively as
$\Phi_{n}$, $\Phi_{n,(p)}$ and $\Phi_{n,p}$
and similarly write $R_n$ and $R_{n,p}$ for
$R_{K_n/k,\Sigma}$ and $R^{(j)}_{K_n/k,\Sigma,p}$.
We write also $V_n$ for the $\bbQ$-vector space
$\bigwedge_{\bbQ G_n}^r\bbQ U_{\Sigma}(K_n)$
and $V_n^0$ for its subspace $V_n^{[\Sigma,G_n,r]}$
which identifies with
$(\bigwedge_{\bbQ G_n}^r\bbQ E(K_n))^{[\Sigma,G_n,r]}$
for $n\geq n_1$, by Lemma~\ref{lemma:1B}.
\emph{Let us assume that $C1(K_n/k,\emptyset)$
holds for all $n\geq n_1$}.
Lemma~\ref{lemma:2A} implies that
it has a unique solution in $V_n^0$,
which we denote $\eta_n$, so that
\beql{eq:2D}
\frac{2^r}{\sqrt{d_k}}R_n(\eta_n)=
\Phi_n(1)\ \ \ \forall n\geq n_1
\eeq
Observe that Lemma~\ref{lemma:1B} also implies that for each $n\geq n_1$ we have
\beql{eq:2D.5}
\mbox{
$\eta_n=\frac{1}{c_n}\bar{\eps}_{1,n}\wedge\ldots\wedge\bar{\eps}_{r,n}$
for some
$c_n\in\bbN$
and
$\eps_{1,n},\ldots\eps_{r,n}\in E(K_n)^{[\Sigma,G_n,r]}$}
\eeq
Now, the sequence $(\eta_n)_{n\geq n_1}$ is
coherent with respect to norms:
Suppose $n\geq m\geq n_1$
and write $N_{n/m}$ and $\pi_{n/m}$ instead of
$N_{K_n/K_m}$ and $\pi_{K_n/K_m}$ respectively.
Then Equation\refeq{1B.5} (with $S=\Sigma$) implies that
$R_m(N_{n/m}\eta_n)=\pi_{n/m}(R_n(\eta_n))$
and Lemma~\ref{lemma:1Z} implies
that $\pi_{n/m}(\Phi_n(1))=\Phi_{m}(1)$.
Applying $\pi_{n/m}$ to\refeq{2D}, it follows that
$N_{n/m}\eta_n$
is a solution of $C1(K_m/k,\emptyset)$.
(Incidentally, this argument also shows
that $C1(K_n/k,\emptyset)$ holds for all  $n\geq n_1$ iff
it holds for infinitely many $n$).
Moreover, if $\chi\in G_m^\dag$ then $\chi\circ\pi_{n/m}\in G_n^\dag$
and it is easy to see that
$r(\Sigma,\chi)=r(\Sigma,\chi\circ\pi_{n/m})$.
From this it follows without difficulty that
$\pi_{n/m}(e_{\Sigma,G_n,>r})=e_{\Sigma,G_m,>r}$.
Thus the $\pi_{n/m}$-semilinearity of $N_{n/m}$ gives
$e_{\Sigma,G_m,>r}N_{n/m}\eta_n=N_{n/m}e_{\Sigma,G_n,>r}\eta_n=0$.
Therefore $N_{n/m}\eta_n$ lies in
$V_n^0$
and by uniqueness it follows that
$
N_{m/n}\eta_n=\eta_m
$.
In particular,
\beql{eq:2Enew}
(\eta_n)_{n\geq
n_1}\in
{\displaystyle\lim_{\longleftarrow}}
\{{\textstyle\bigwedge^r_{\bbQ G_n}}\bbQ
E(K_n):n\geq n_1\}
\subset
{\displaystyle\lim_{\longleftarrow}}\{{\textstyle\bigwedge^r_{\bbQ G_n}}\bbQ
U_\Sigma(K_n):n\geq n_1\}
\eeq
where the limits are taken with respect to the norm maps $N_{n/m}$.

The core idea of this paper is to
investigate the consequences of a stronger hypothesis than\refeq{2Enew}
in which, very roughly speaking, one reverses the
order of the functors $\displaystyle\lim_{\longleftarrow}$ and
$\bigwedge^r\bbQ\otimes_\bbZ$ on the right hand side.
Now, on the one hand
$\displaystyle\lim_{\longleftarrow}U_{\Sigma}(K_n)$
equals $\displaystyle\lim_{\longleftarrow}U_{T_p}(K_n)$
(since the primes above any finite prime
$\fq\nin T_p$ are inert and unramified in
$K_\infty/K_N$ for some $N=N(\fq)$).
On the other, the consequences we have in mind will
require norm-coherent
sequences of elements which are units \emph{locally}
at primes above $T_p$. We therefore consider the existence of elements of
$\bigwedge_{\bbQ G_\infty}^r\bbQ \cE_\infty(K)$
where $\cE_\infty(K):=\displaystyle\lim_{\longleftarrow}E(K_n)$
and $\bbQ \cE_\infty (K)=\bbQ\otimes_\bbZ \cE_\infty (K)$ is treated as
a module for the (uncompleted) group-ring $\bbQ G_\infty$.
We shall write  $\bar{\ueps}$ for the image of an element
$\ueps\in \cE_\infty(K)$ in $\bbQ\cE_\infty(K)$ and,
for each $n\geq n_1$, we denote by $\beta_n$ the homomorphism
from $\bigwedge_{\bbQ \Ginf}^r\bbQ \cE_\infty (K)$ to
$\bigwedge_{\bbQ G_n}^r\bbQ E(K_n)$ (both notated additively) which sends
$\bar{\ueps}_1\wedge\ldots\wedge\bar{\ueps}_r$ to
$\bar{\eps}_{1,n}\wedge\ldots\wedge\bar{\eps}_{r,n}$.
For any $K/k$ and $p$ satisfying\refeq{2A},\refeq{2B} and\refeq{2C}
we formulate the
\begin{property}{$P1(K/k,p)$}
There exists $\ueta\in\bigwedge_{\bbQ \Ginf}^r\bbQ\cE_\infty(K)$
such that, for every $n\geq n_1$, $\beta_n(\ueta)$ lies in $V_n^0$
and is a solution of
$C1(K_n/k,\emptyset)$ (the unique solution in $V_n^0$) \ie\
\[
\Phi_n(1)=\frac{2^r}{\sqrt{d_k}}
R_n(\beta_n(\ueta))\ \ \ \forall n\geq n_1
\]
(We shall say that $\ueta$
\emph{demonstrates} $P1(K/k,p)$).
\end{property}
\noindent\rem\ It is unclear at present whether we should expect this for every
such pair $(K/k,p)$, which is why we have called
it a `property' (that may or may not be possessed) rather than a conjecture.\vspace{1ex}\\
\rem\label{rem:2B}\
Any $\ueta$ demonstrating $P1(K/k,p)$ may be written as a finite sum
\beql{eq:2E.25}
\ueta=\sum_{w=1}^W\frac{1}{c_w}\bar{\ueps}^{(w)}_1
\wedge\ldots\wedge\bar{\ueps}^{(w)}_r
\eeq
where each $c_w$ is a positive integer and each
$\ueps^{(w)}_i$ is a norm-coherent sequence $(\eps^{(w)}_{i,n})_{n\geq n_1}$
with $\eps^{(w)}_{i,n}\in E(K_n)$ for each $n$.
The condition $\beta_n(\ueta)\in V_n^0$
implies that $\eta_n=\beta_n(\ueta)$ satisfies
\refeq{2D.5} for each $n\geq n_1$. However,
the same assumption at the infinite level
(namely that we may take $W=1$ in\refeq{2E.25}) appears to be a strict
strengthening of $P1$. We shall not consider it further here
except to note that it does hold
in all the cases with $K/\bbQ$ abelian
so far considered (see Example~\ref{ex:2B} below).
\vspace{1ex}\\
\rem\label{rem:denoms}\ Let us define the \emph{denominator}
of an element $\eta\in\bigwedge^r_{\bbQ G_n}\bbQ E(K_n)$ (for any $n$) to be the
smallest positive integer $d$ such that $d\eta$
lies in the lattice which is the image of $\bigwedge^r_{\bbZ G_n}E(K_n)$
in $\bigwedge^r_{\bbQ G_n}\bbQ E(K_n)$. Property $P1(K/k,p)$ implies in particular
that the unique solutions of $C1(K_n/k,\emptyset)$ lying in $V^0_n$ have \emph{bounded
denominators} as $n$ increases. (Indeed, if $\ueta$ is given by\refeq{2E.25}
and demonstrates $P1(K/k,p)$, then for each
$n\geq n_1$,
the denominator of $\beta_n(\ueta)$ clearly divides the l.c.m.\ of the $c_w$'s.)
\vspace{1ex}\\
We consider some special cases of $P1(K/k,p)$.\vspace{1ex}\\
\example\label{ex:2A}\ {\bf The Case $k=\bbQ$.}\ In this case $r=1$,
$K$ is an abelian field and
Condition\refeq{2B} implies that
$\ff(K)=f'p^{n_0}\bbZ$ for some $f'\in\bbZ_{>1}$ prime to $p$
and $n_0\geq 0$. For every $n\geq n_1$ we have $\ff(K_n)=f_n\bbZ$ where
$f_n:=f'p^{\max(n_0,n+1)}$ so that $k(\ff(K_n))$
is the real cyclotomic field $\bbQ(\zeta_{f_n})^+$.
For each $n\geq n_1$ we let
$\eps_n=\eps_{1,n}$ be the norm from $\bbQ(\zeta_{f_n})^+$ to $K_n$ of
$(1-\zeta_{f_n})(1-\zeta_{f_n}\inv)$.
Then it can be shown that $\eta_n:=-\frac{1}{2}\bar{\eps}_n$
lies in $V_n^0=(\bbQ E(K_n))^{[\Sigma,G_n,1]}$ and
that if we take $\ttau=\ttau_1=1\in\Gal(\barbbQ/\bbQ)$ then $\eta_n$ is a solution of
$C1(K_n/\bbQ,\emptyset)$. (See \eg~\cite[\S~3.5]{zetap1}
for the case in which $K_n$ equals $\bbQ(\zeta_{f_n})^+$. The general case follows
by `taking norms'.) Now, is easy to check that the sequence
$(\eps_n)_{n\geq n_1}$ is norm-coherent (in fact, this coherence is
already a consequence of\refeq{2Enew}, at least up to $\pm1$). Therefore
$\ueta:=-\frac{1}{2}\otimes(\eps_n)_{n\geq n_1}$ lies in $\bbQ\cE_\infty(K)$ and
$\eta_n=\beta_n(\ueta)$. So $\ueta$
demonstrates $P1(K/\bbQ,p)$.\vspace{1ex}\\
\example\label{ex:2B}\ {\bf The Case $K/\bbQ$ Abelian.} This
generalises the previous case inasmuch as it does not assume that
$k=\bbQ$. In this case, we shall use the methods of~\cite[\S 3]{Po2}
to sketch a proof that
$P1(K/k,p)$ holds under\refeq{2A},\refeq{2B},\refeq{2C}
\emph{together with
the additional hypothesis}
\beql{eq:2E.5}
\mbox{$\Sigma(K/k,p)$ is precisely the set of places in $k$ above
$\Sigma(K/\bbQ,p)$}
\eeq
This allows us to `base change' from the previous example, using induction.
Note first that the fixed, extended embeddings
$\ttau_i$ for  $i=1,\ldots,r$
(used to define the $\lambda_{K_n/k,i}$ for any $n$)
form a set of coset representatives
for $\Gal(\barbbQ/k)$ in $\Gal(\barbbQ/\bbQ)$. For each $n\geq n_1(K,p)$ the group-ring
$\bbC\Gal(K_n/\bbQ)$ may be considered as a free module
over $\bbC G_n=\bbC \Gal(K_n/k)$
with basis provided by $\{\ttau_1|_{K_n},\ldots,\ttau_r|_{K_n}\}$.
For any $x\in \bbC\Gal(K_n/\bbQ)$
we write $\det_{\bbC G_n}(x)\in\bbC G_n$ for the determinant of
multiplication by $x$ and, to avoid confusion, we shall
denote by
$\eta_n(K/\bbQ)$ and $\ueta(K/\bbQ)$ respectively
the elements $\eta_n\in\bbQ E(K_n)$ and $\ueta\in\bbQ\cE_\infty(K)$
described in the previous example.
Using the fact that $\eta_n(K/\bbQ)$ is a solution of $C1(K_n/\bbQ,\emptyset)$ and
calculating $\det_{\bbC G_n}$ by means of the basis
described above, we find easily that
\begin{eqnarray}
{\textstyle\det_{\bbC G_n}(\half}\Phi_{K_n/\bbQ,\emptyset}(1))&=&
{\textstyle\det_{\bbC G_n}}(\lambda_{K_n/\bbQ,1}(\eta_n(K/\bbQ)))\nonumber\\
&=&\det(\lambda_{K_n/k,i}(\ttau_t\inv\eta_n(K/\bbQ)))_{i,t=1}^r\nonumber\\
&=&R_{K_n/k}
(\ttau_1\inv\eta_n(K/\bbQ)\wedge\ldots\wedge\ttau_r\inv\eta_n(K/\bbQ))
\label{eq:2E.75}
\end{eqnarray}
for all $n\geq n_1$ (provided that $\lambda_{K_n/\bbQ,1}$ is defined
taking `$\ttau_1$ for $\bbQ$' to be $1$, as in the previous example).
On the other hand, it follows from~\cite[Thm. 5.1]{twizas} that for $n\geq n_1$, the elements
$\half\Phi_{K_n/\bbQ,\emptyset}(1)$ and $\frac{\sqrt{d_k}}{2^r}\Phi_{K_n/k,\emptyset}(1)$
are equal respectively to the elements which would be denoted
$\Theta'_{K_n/\bbQ,S}(0)$ and $\Theta^{(r)}_{K_n/k,S'}(0)$ in~\cite{Po2},
where $S=\Sigma(K/\bbQ,p)$ and $S'=\Sigma(K/k,p)$.
Condition\refeq{2E.5} therefore puts us in the situation considered by Popescu and his
Equation~(3) implies that
$\Theta^{(r)}_{K_n/k,S'}(0)=\det_{\bbC G_n}(\Theta'_{K_n/\bbQ,S}(0))$. Putting this
together with\refeq{2E.75} shows that
\[
\frac{\sqrt{d_k}}{2^r}\Phi_{K_n/k,\emptyset}(1)=
R_{K_n/k}
(\ttau_1\inv\eta_n(K/\bbQ)\wedge\ldots\wedge\ttau_r\inv\eta_n(K/\bbQ))
\]
Also, since $\eta_n(K/\bbQ)$ lies in
$(\bbQ E(K_n))^{[\Sigma(K_n/\bbQ,p),\Gal(K_n/\bbQ),1]}$
and\refeq{2E.5} holds, an argument in~\cite{Po2} shows
that $\ttau_1\inv\eta_n(K/\bbQ)\wedge\ldots\wedge\ttau_r\inv\eta_n(K/\bbQ)$ lies
in $(\bigwedge_{\bbQ G_n}^r\bbQ E(K_n))^{[\Sigma(K_n/k,p),G_n,r]}$. Hence
it is a solution of $C1(K_n/k,\emptyset)$ whenever $n\geq n_1$ and it follows immediately
that $\ttau_1\inv\ueta(K/\bbQ)\wedge\ldots\wedge\ttau_r\inv\ueta(K/\bbQ)$
is an element of $\bigwedge_{\bbQ G_\infty}^r\bbQ \cE_\infty(K)$ demonstrating $P1(K/k,p)$.\vspace{2ex}\\
We have also verified $P1(K/k,p)$
in certain cases with $K/\bbQ$ abelian for which
which\refeq{2E.5} does not hold and
it is quite possible that the latter hypothesis
is in actually unnecessary. On the other hand, the author
knows of no extension of real number fields $K/k$
with $K\not\subset\bbQ^{\rm ab}$ and
such that even the basic Stark Conjecture
$C1(K_n/k,\emptyset)$ has been proven
for infinitely many $n$, for some $p$. This makes it hard to
construct further examples in which $P1$ can be shown either to hold
or to fail.

Our main object now is to study a \emph{$p$-adic} analogue of
property $P1(K/k,p)$, in which, roughly
speaking, we replace $C1(K_n/k,\emptyset)$ by $C2(K_n/k,T_p,p)$.
We must therefore start by examining the passage from
$T=\emptyset$ to $T=T_p$ for the complex conjecture $C1$
and this in turn requires a
closer analysis of ramification above $p$ in the extension $K_\infty/k$.
Let $L/F$ be a finite abelian field
extension. If $L$ and $F$ are local fields and $v\in[-1,\infty)$,
we recall that $G(L/F)^v$ denotes the $v$th ramification group
(in the `upper numbering', see~\cite[Ch.~IV]{SerreLF}).
If $L$ and $F$ are number fields and $\fp$ is a prime of $\cO_F$,
then $G(L/F)^v_\fp$ denotes the $v$th ramification group
at any prime $\fP$ of $\cO_L$ above $\fp$, naturally identified
with the group $G(L_\fP/F_\fp)^v$ of the completed
extension.
The upper numbering is preserved under restriction
and (since $L/F$ is abelian)
depends only on $\lceil v \rceil$, the `integer ceiling' of $v$.
Taking $F$ to be $k$, local and global class field theory
give the formula
\beql{eq:2F}
\ord_\fq(\ff(L))=
\min\{v\in\bbN:G(L/k)^v_\fq=\{1\}\}\ \ \ \mbox{for all primes $\fq$ of $\cO$}
\eeq
\begin{lemma}\label{lemma:2B} Suppose that $n,v\geq 1$ are integers
and let $\fp|p$.
If $v\leq n$ then $G(k_n/k)^v_\fp=\Gal(k_n/k_{v-1})$
which is cyclic of order $p^{n-v+1}$. Otherwise,
$G(k_n/k)^v_\fp=\{1\}$.
\end{lemma}
\bPf\
Let $\widehat{k_n}$ and $\widehat{\bbB_n}$ denote the completions
of $k_n$ and $\bbB_n$ at the unique primes above $\fp$ and $p$
respectively. Then $\widehat{k_n}$ is the compositum
of $\widehat{\bbB_n}$ with $k_\fp$. But
$k_\fp/\bbQ_p$ is unramified by\refeq{2C}. We deduce that
$\widehat{k_n}/\bbQ_p$ is abelian and, easily, that the group
$G(\widehat{k_n}/k_\fp)^v$ equals $G(\widehat{k_n}/\bbQ_p)^v$.
Thus it is a subgroup of $\Gal(\widehat{k_n}/k_\fp)$
mapped isomorphically to
$G(\widehat{\bbB_n}/\bbQ_p)^v$ by restriction to
$\widehat{\bbB_n}$.
But the latter group is also the image
of $G(\bbQ_p(\zeta_{\pnpo})/\bbQ_p)^v$ which equals
$\Gal(\bbQ_p(\zeta_{\pnpo})/\bbQ_p(\zeta_{p^v}))$
if $v\leq n$, $\{1\}$ if not
(\cite[Ch.~IV]{SerreLF}). Thus if $v\leq n$ then
$G(\widehat{\bbB_n}/\bbQ_p)^v=
\Gal(\widehat{\bbB_n}/\widehat{\bbB_{v-1}})$ and so
$G(\widehat{k_n}/k_\fp)^v$ equals
$\Gal(\widehat{k_n}/\widehat{k_{v-1}})$ and is isomorphic to
$\bbZ/p^{n-v+1}\bbZ$. Otherwise, $G(\widehat{k_n}/k_\fp)^v=\{1\}$.
Of course, $G(\widehat{k_n}/k_\fp)^v$ and
$\Gal(\widehat{k_n}/\widehat{k_{v-1}})$ are the images of
$G(k_n/k)_\fp^v$ and
$\Gal(k_n/k_{v-1})$ respectively under the natural
identification of
$\Gal(k_n/k)$ with $\Gal(\widehat{k_n}/k_\fp)$, so the result follows. \ePf
\noindent Since $k_n/k$ is unramified
outside $p$, Lemma~\ref{lemma:2B},\refeq{2F} and\refeq{2C} give
$\ff(k_n)
=p^{n+1}\cO$ for all $n\geq 1$. If
we define $\ff_n$ to be $\ff(K_n)$ for all $n\geq 0$
(so $\ff_0=\ff(K)$) then it follows easily that
\beql{eq:2G}
\ff_n={\rm l.c.m.}(\ff_0,\ff(k_n))={\rm l.c.m.}(\ff_0,p^{n+1}\cO)\ \ \
\mbox{for all $n\geq 1$}
\eeq
Now define
\[
n_2=n_2(K)=n_2(K/k,p)=\max(\{1\}\cup\{\ord_\fp(\ff_0):\fp|p\})=\max(\{n_1\}\cup\{\ord_\fp(\ff_0):\fp|p\})
\]
(For the fourth equality, note that $n_1\leq 1$ and
that $\ord_\fp(\ff_0)=0\ \forall\ \fp|p$ implies
that $K/k$ is unramified above
$p$ and hence that $n_1=1$.)
Equation\refeq{2G} implies that
\beql{eq:xxtra}
\mbox{the $T_p$-part
of $\ff_n$ is $p^{n+1}\cO$ for every $n\geq \max(1, n_2-1)$}
\eeq
(In particular
for every $n\geq n_2$).\vspace{1ex}\\
\rem\ We are \emph{not} assuming that the extensions $K/k$ and $\Kinf/k$
are linearly disjoint. Thus if we define $n_3\geq 0$ by $K\cap\kinf=k_{n_3}$,
then it may be that $n_3\geq 1$, in which case $p^{n_3+1}\cO=\ff(k_{n_3})$
divides $\ff_0$. Thus in all cases $n_2\geq n_3+1$, although this inequality
may very well be strict.
\begin{lemma}\label{lemma:2C}
Suppose that $n,v\geq 1$ are integers with $v\geq n_2$ and let $\fp|p$.
If $v\leq n$ then
$G(K_n/k)_\fp^v=\Gal(K_n/K_{v-1})$ which is cyclic of order $p^{n-v+1}$.
Otherwise $G(K_n/k)_\fp^v=\{1\}$.
\end{lemma}
\bPf\
Equation\refeq{2F} and
the definition of $n_2$ imply that
$G(K/k)^v_\fp=1$. Hence, restricting
$G(K_n/k)_\fp^v$ to $K$ we find that it lies in $\Gal(K_n/K)$. The
restriction to $k_n$ therefore maps it isomorphically
onto $G(k_n/k)^v_\fp$. Now apply Lemma~\ref{lemma:2B}.
\ePf\noindent
For any $\chi\in G_n^\dag$ and any prime $\fq$ of $\cO$
we clearly have
$\ord_\fq(\ff(\chi))=\min\{v\in\bbN:G(K_n/k)^v_\fq\subset\ker \chi\}$.
Taking $v=n$ in the above lemma, we see that \emph{if $n\geq n_2$} then
$[K_n:K_{n-1}]=p$ and, for every $\chi\in G_n^\dag$,
\beql{eq:2H}
(\ff(\chi)|\fp\inv\ff_n\ \mbox{for some $\fp|p$})
\Longleftrightarrow
\chi(\Gal(K_n/K_{n-1}))=1
\Longleftrightarrow
\ff(\chi)|p
\inv\ff_n
\eeq
For each $n\geq n_2$
we write $e_n$ for the idempotent of $\bbQ G_n$
given by $1-p\inv\sum_{\sigma\in\Gal(K_n/K_{n-1})}\sigma$
and note that $\ff_n$ is not
$T_p$-trivial since $p^2|\ff_n$ (or by\refeq{2B}).
Thus $\Phi_{n,(p)}:=\Phi_{K_n/k,T_p}$
is holomorphic.
\begin{prop}\label{prop:2A}
If $n\geq n_2$ then $\Phi_{n,(p)}(s)=
p^{r(n+1)(1-s)}e_n\Phi_n(s)$
for all $s\in \bbC$.
\end{prop}
\bPf\
It suffices to prove that
\beql{eq:2I}
\chi(\Phi_{n,(p)}(s))=p^{r(n+1)(1-s)}\chi(e_n\Phi_n(s))
\ \ \ \forall\chi\in G_n^\dag,\ \forall s\in\bbC\setminus\{1\}
\eeq
We use Proposition~\ref{prop:1A}
with $K=K_n$ (so $\fh=\ff_n\ff(\chi)\inv$)and suppose first that there
exists $\fp|p$ such that $\ff(\chi)|\fp\inv\ff_n$.
The case $T=T_p$ of the said Proposition then has
$\ord_\fp(\fh)\geq 1$ but also $\ord_\fp(\ff_n)\geq 2$ from which it follows
that $\fh_0$ is either not square-free or not prime to $\ff(\chi)$, hence the
L.H.S. of\refeq{2I} equation vanishes (for all $s\in\bbC$).
But\refeq{2H}
implies that $\chi(e_n)=0$ so the R.H.S. vanishes as well.
If no such $\fp$ exists
then the $T_p$-part of $\ff(\chi)$ is $p^{n+1}\cO$ and in both the
cases $T=T_p$ and $T=\emptyset$ of Proposition~\ref{prop:1A}
every $\fp\in T$ divides $\ff(\chi)$ but neither $\fh$ nor (by\refeq{2C})
$\cD$. Therefore, we can apply\refeq{1G}
in both cases to get (for $s\neq 1$):
\begin{eqnarray*}
\chi(\Phi_{n,(p)}(s))&=&
g(\hchi)p^{r(n+1)(1-s)}N(\ff(\chi)\cD)^{s-1}
\prod_{\fp\ \mathit{prime}\atop \fp|\fh,\ \fp\ndiv\ff(\chi)}
\left(1-\frac{\hchi\inv(\fp)}{N\fp^{1-s}}\right)L(s,\hchi)\\
&=&p^{r(n+1)(1-s)}\chi(\Phi_n(s))
\end{eqnarray*}
But\refeq{2H} now implies that $\chi(e_n)=1$ in this case.
So, again,\refeq{2I}
holds.\ePf
\noindent
If $\ueta$ demonstrates $P1(K/k,p)$
then for each $n\geq n_2$, Proposition~\ref{prop:2A} gives
$\Phi_{n,(p)}(1)=e_n\Phi_n(1)=
(2^r/\sqrt{d_k})R_n(e_n\beta_n(\ueta))$
and since\refeq{2C} implies
$\prod_{\fp\in T_p}N\fp=p^r$, it follows that,
for each $n\geq n_2$,
$p^re_n\beta_n(\ueta)$
is a solution of the complex conjecture $C1(K_n/k,T_p)$ lying
in $V_n^0$.
To deduce that it is also a solution
of the $p$-adic conjecture $C2(K_n/k,T_p,p)$ it is necessary
and sufficient to assume only the existence of \emph{some} common
solution $\eta_n'$
of $C1(K_n/k,T_p)$ and $C2(K_n/k,T_p,p)$ lying in $V_n$,
that is, a solution of $C3(K_n/k,T_p,p)$.
To see this, we use the following claim (proved below)
\begin{claim}\label{claim:2A}
For any such $\eta_n'$ (with $n\geq 1$) the element $e_{\Sigma,G_n,r}\eta_n'$
is also a common solution of $C1(K_n/k,T_p)$ and $C2(K_n/k,T_p,p)$.
\end{claim}
Indeed, assuming this, $p^re_n\beta_n(\ueta))$ must be \emph{equal} to
$e_{\Sigma,G_n,r}\eta'_n$
(since both lie in $V_n^0$
on which $R_n$ is injective), so the former element is also
solution of $C2(K_n/k,T_p,p)$, as required. To sum up, we formulate
a $p$-adic property potentially satisfied
by any $K/k$ and $p$ satisfying\refeq{2A},\refeq{2B} and\refeq{2C}.
\begin{property}{$P2(K/k,p)$}
There exists $\ueta\in\bigwedge_{\bbQ \Ginf}^r\bbQ\cE_\infty(K)$
such that, for every $n\geq n_2$, $\beta_n(\ueta)$ lies in $V_n^0$
and $p^re_n\beta_n(\ueta))$ is a solution of
$C2(K_n/k,T_p,p)$, \ie\
\[
\Phi_{n,p}(1)=\frac{2^r}{j(\sqrt{d_k})}
R_{n,p}(e_n\beta_n(\ueta))\ \ \ \forall n\geq n_2
\]
(We shall say that $\ueta$
\emph{demonstrates} $P2(K/k,p)$).
\end{property}
\noindent
The preceding argument gives
\begin{prop}\label{prop:2B}
Assume that $C3(K_n/k,T_p,p)$ holds for all $n\geq n_2$. If
$\ueta$ demonstrates $P1(K/k,p)$
then it also demonstrates $P2(K/k,p)$.
\end{prop}
\bPf\ It remains only to give the \emph{Proof of Claim~\ref{claim:2A}:}
It suffices to show
(\cf\ the proof of Lemma~\ref{lemma:2A}) that
$\chi(\Phi_{n,(p)}(1))=\chi(\Phi_{n,p}(1))=0$ for any $\chi\in G_n^\dag$
with $r(\Sigma,\chi)>r$. But for any such $\chi$ the set
$
\Sigma(\chi):=\{\fp\in\Sigma\,:\,\chi(G(\fp))=1\}
$
is non-empty. Moreover, all its elements are clearly prime to $\ff(\chi)$.
Now apply Proposition~\ref{prop:1A} and Corollary~\ref{cor:1A}
with $K=K_n$, $T=T_p$.
Suppose first that there exists $\fp\in\Sigma(\chi)\cap T_p$
(for example, if $\chi=\chi_0$).
In this case,
$\ord_\fp(\fh_0)=\ord_\fp(\fh)=\ord_\fp(\ff(K_n))\geq 2$
by\refeq{2G}
and the hypothesis that $n\geq 1$, so that $\chi(\Phi_{n,(p)}(s))$
and $\chi(\Phi_{n,p}(s))$ vanish identically.
If, on the other hand, $\Sigma(\chi)\cap T_p=\emptyset$ then
$\chi\neq\chi_0$ and each $\fp\in\Sigma(\chi)$ divides
$\fh'$ thus contributing a factor vanishing at $s=1$
to the first product on the right-hand sides of\refeq{1G} and\refeq{1H}.
The Claim follows.
\ePf
\rem\ If $\ueta$ is of the form given in\refeq{2E.25}
then it is easy to see that for each $n\geq n_2$,
\[
p^re_n\beta_n(\ueta)=\sum_{w=1}^W\frac{1}{c_w}\bar{\eps'}^{(w)}_{1,n}
\wedge\ldots\wedge\bar{\eps'}^{(w)}_{r,n}
\]
where $\eps'^{(w)}_{i,n}$ denotes
$(\eps^{(w)}_{i,n})^p/N_{n/(n-1)}\eps^{(w)}_{i,n}=
(\eps^{(w)}_{i,n})^p/\eps^{(w)}_{i,n-1}$ (the last equality only makes sense
if $n-1\geq n_1$).
\section{Consequences of Property $P2$}\label{section:consP2}\
We shall study these by means of the $p$-adic power series
attached by Coleman~\cite{Coleman} to norm-coherent sequences in
certain towers of local fields. First, some basic facts and notations. Let
$n\in\bbZ$, $n\geq -1$. Since $p$ is fixed, we shall usually abbreviate
$\zeta_\pnpo$ to $\zeta_n$ (whenever this abuse cannot cause confusion)
and identify this element of $\barbbQ$
with its image in $\bbC_p$ under $j$. For any subfield $L$ of
$\bbC$ or $\bbC_p$ we shall write $\tilde{L}_n$ for $L(\zeta_n)$
(so $\tilde{L}_{-1}=L$)
and $\tilde{L}_\infty$ for $\bigcup_{n\geq-1} \tilde{L}_n$.
For any closed subfield $L$ of $\bbC_p$ we shall
write $\cA(L)$ for the
subspace of $L[[X]]$ consisting of those power series
$h(X)=\sum_{i=0}^\infty a_iX^i$
whose norm $||h||:=\sup\{|a_i|_p:i\geq 0\}$ is finite.
It is is an $L$-Banach space under $||\,\cdot\,||$.
For any closed and open subset $U$ of $\bbZ_p$,
we denote by ${\rm Meas}(U,L)$ the $L$-Banach space of of all (bounded)
$L$-valued measures on $U$ equipped with the natural norm.
`Extension by zero' allows us to regard ${\rm Meas}(U,L)$ as the subspace of
${\rm Meas}(\bbZ_p,L)$ consisting of those measures
supported on $U$. Moreover, there is a well known correspondence
between $\cA(L)$ and ${\rm Meas}(\bbZ_p,L)$ which
is in fact an isometry of Banach spaces.
(More details of the correspondence may be found
in~\cite[App. 5,6]{Sch}, \cite{La} or~\cite{alshig2}.
The $p$-adic `Weierstrass Preparation Theorem' implies the
following well known and useful fact
(see \eg~\cite[Cor. 7.4]{Wa}):\vspace{1ex}\\
{\bf Uniqueness Principle}\ \
If $[L:\bbQ_p]$ is finite then a non-zero
element of $\cA(L)$ can have only finitely many roots
in $\{a\in\bbC_p:|a|_p<1\}$.\vspace{1ex}\\
Let $H$ be a finite, unramified extension of $\bbQ_p$ so that
$\{\tilde{H}_n\}_{n\geq -1}$ is the division tower over $H$
associated to the multiplicative
Lubin-Tate formal group (with group law $X+Y+XY$) over $\bbQ_p$.
The element $\pi_n:=\zeta_n-1$ is a uniformiser of $\tilde{H}_n$
for all $n\geq 0$ and $[p]\pi_n=(1+\pi_n)^p-1=\pi_{n-1}$ for
$n\geq 1$.
The Frobenius element $\phi$ of $\Gal(H/\bbQ_p)$
extends uniquely to a continuous automorphism of $\tH_\infty$
fixing $\zeta_n$ for all $n$ and
acts coefficientwise on
the ring $\tH_\infty[[X]]$ of formal power series.
We let $\hat{\kappa}$ denote the
Lubin-Tate isomorphism from $\Gal(\tH_\infty/H)$ to
$\bbZ_p^\times$ so that $\sigma(\pi_n)=[\hat{\kappa}(\sigma)]\pi_n=
(1+\pi_n)^{\hat{\kappa}(\sigma)}-1$ for all $\sigma\in\Gal(\tH_\infty/H)$ and all $n$.
Thus $\Gal(\tH_\infty/H)$ also acts on $\bbC_p[[X]]$ by letting
$\sigma\cdot f(X):=f((1+X)^{\hat{\kappa}(\sigma)}-1)$.
Let us write $U_\infty(H)$ for
${\displaystyle\lim_{\longleftarrow}\{U(\tH_n):n\geq 0\}}$, the
inverse limit of the local units with respect to the
norm maps $N_{\tH_m/\tH_n}$. For each
$\uu=(u_n)_{n\geq 0}$ in $U_\infty(H)$, it follows
from~\cite[Thm. A]{Coleman}
that there exists a formal power series $g\in\cO_H[[X]]^\times$ such that
$\phi^{-n}g(\pi_n)=u_n$ for all $n\geq 0$ and by the Uniqueness Principle,
any infinite subsequence of these equations determines
$g$ in $\cO_H[[X]]$. Together with norm-coherence this
implies that $g$ must satisfy the relation
\beql{eq:3A} \prod_{\zeta\in\mu_p}g(\zeta(1+X)-1)=\phi g((1+X)^p-1)\eeq
and be independent of the choice of
unramified $H$ such that $\uu\in U_\infty(H)$.
We therefore denote it simply $g(\uu;X)$. Moreover $\Gal(\tH_\infty/H)$
acts naturally on $U_\infty(H)$ and the Uniqueness Principle also
implies that
\beql{eq:3A.5}
g(\sigma\uu;X)=\sigma\cdot g(\uu;X)\ \ \
\forall\uu\in U_\infty(H),\ \forall \sigma\in\Gal(\tH_\infty/H)
\eeq
Let
$g^\ast(\uu;X)$ denote the power series
$g(\uu;X)^p/\phi g(\uu;(1+X)^p-1)$. Simple arguments show that it satisfies
$\prod_{\zeta\in\mu_p}g^\ast(\uu;\zeta(1+X)-1)=1$ and lies in $1+p\cO_H[[X]]$.
Therefore,
\[
h^\ast(\uu;X):=\frac{1}{p}\log(g^\ast(\uu;X))=\sum_{i=1}^\infty(-1)^{i-1}(g^\ast(\uu;X)-1)^i/pi
\]
defines an element of $\cO_H[[X]]$. It follows from the
above definitions that
$\phi^{-n}h^\ast(\uu;\pi_n)=\frac{1}{p}\log_p(u_n^\ast)$ for all $n\geq 1$
where $u_n^\ast:=u_n^p/u_{n-1}$, and also
that $h(X)=h^\ast(\uu;X)$ satisfies the relation
\beql{eq:3B}
\sum_{\zeta\in\mu_p}h(\zeta(1+X)-1)=0
\eeq
If $L$ is any closed subfield of $\bbC_p$, we let
$\cA^\ast(L)$ denote the subspace of $\cA(L)$ consisting of all
power series $h\in\cA(L)$ satisfying\refeq{3B}. Under the
power series/measure correspondence mentioned above,
$\cA^\ast(L)$ corresponds to ${\rm Meas}(\bbZ_p^\times,L)$.
Two maps from ${\rm Meas}(\bbZ_p,L)$ to itself are defined
by the operation of multiplying a measure respectively by
the functions $f_1:x\mapsto x$ and $f_2:x\mapsto\omega\inv(x)$
(the latter extended by zero from $\bbZ_p^\times$ to $\bbZ_p$).
It is easy to check that both maps restrict to isometries from
${\rm Meas}(\bbZ_p^\times,L)$ to itself.
Multiplication by $f_1$ corresponds to the differential operator
$D:h(X)\mapsto (1+X)h'(X)$ on $\cA(L)$,
while multiplication by $f_2$ is easily seen to correspond to the map
from $\cA(L)$ to itself given by
\[
\mbox{$V:h(X)\longmapsto{\displaystyle
-\frac{g(\omega)}{p}
\sum_{l=1}^{p-1}}\omega(l)h(\zeta_0^l(1+X)-1)$
\ \ \ where\ \ \
$g(\omega):={\displaystyle\sum_{l=1}^{p-1}}\zeta_0^l\omega\inv(l)$}
\]
It follows that $D$ and $V$ restrict to commuting
isometries from $\cA^\ast(L)$ to itself (under the $||\,\cdot\,||$-norm)
as does their composite $\check{D}=D\circ V=V\circ D:\cA(L)\rightarrow\cA(L)$,
corresponding as it does to the
multiplication of a measure by the function
$f_1f_2:x\mapsto\langle x\rangle$, extended by zero.
For any $h\in \cA^\ast(L)$ and
any $t\in\bbZ$, positive or negative, we shall therefore
write simply $D^{t}h$ ({\em resp.} $\check{D}^{t}h$)
for $(D|_{\cA^\ast(L)})^{t}h$ ({\em resp.} $(\check{D}|_{\cA^\ast(L)})^{t}h$)
which lies in $\cA^\ast(L)$ but is obviously
independent of the choice of $L$ given $h$. Clearly,
$D^{t}h=\check{D}^{t}h$ whenever $(p-1)|t$. We can now define
some twisted and generalised versions of the Coates-Wiles homomorphisms
(compare~\cite[Ch. 13]{Wa} for example).

If $\uu$ lies in
$U_\infty(H)$ then $\check{D}^t h^\ast(\uu;X)$ lies in $\cO_H[[X]]$.
For each $t,n\in\bbZ$ with $n\geq 1$ we may therefore define a map
\displaymapdef{\check{\delta}_{t,n}}{U_\infty(H)}{\cO_{\tH_n}}{\uu}
{\phi^{-n}\check{D}^t h^\ast(\uu;\pi_n)}
It is easy to see that $\check{\delta}_{t,n}$ is a homomorphism
(\ie\ $\check{\delta}_{t,n}(\uu_1\uu_2)=
\check{\delta}_{t,n}(\uu_1)+\check{\delta}_{t,n}(\uu_2)$)
and that
\beql{eq:3B.5}
\check{\delta}_{0,n}(\uu)=\frac{1}{p}\log_p(u_n^\ast)\ \forall\,n\geq 1
\eeq
However, for $t\neq 0$, the value of $\check{\delta}_{t,n}(\uu)$ depends
on the entire sequence $(u_i^\ast)_{i\geq 1}$ (or $(u_i)_{i\geq 0}$), not just
on $u_n^\ast$, indeed we have:
\begin{lemma}\label{lemma:3Z}
Let $\sigma\in\Gal(\tH_\infty/H)$,
$\uu\in U_\infty(H)$ and suppose that $t,n\in\bbZ$ with $n\geq 1$.
Then $\check{D}^t h^\ast(\sigma\uu;\pi_n)=\langle\hat{\kappa}(\sigma)\rangle^t
\check{D}^t h^\ast(\uu;\sigma(\pi_n))$. In particular, if
$\sigma$ lies in $\Gal(\tH_\infty/\tH_n)$ then
$\check{\delta}_{t,n}(\sigma\uu)=\langle\hat{\kappa}(\sigma)
\rangle^{t}\check{\delta}_{t,n}(\uu)$.
\end{lemma}
\bPf\ For each $\sigma\in\Gal(\tHinf/H)$ and
$h\in\cA^\ast(H)$ one checks that $D\sigma\cdot h=\hat{\kappa}(\sigma)\sigma\cdot Dh$
and $V\sigma\cdot h=\omega\inv(\hat{\kappa}(\sigma))\sigma\cdot
Vh$ so that $\check{D}\sigma\cdot h=
\langle\hat{\kappa}(\sigma)\rangle\sigma\cdot \check{D}h$.
It follows that $\check{D}^t\sigma\cdot h=
\langle\hat{\kappa}(\sigma)\rangle^t\sigma\cdot \check{D}^t h$
for all $t\in\bbN$. Since $\sigma$ and $\check{D}$ map
$\cA^\ast(H)$ bijectively to itself,
this equation also holds for $t\in\bbZ$.
Equation\refeq{3A.5}
implies that $h^\ast(\sigma\uu;X)=\sigma\cdot h^\ast(\uu;X)$ and the
first statement follows. The second statement clearly follows from the first by applying
$\phi^{-n}$.\ePf

We now return to the global situation and the notations of the last section,
with $K/k$ and $p$ satisfying\refeq{2A} and\refeq{2B}. For the
purposes of the present section however, we replace\refeq{2C} by
the stronger hypothesis
\beql{eq:3C} \mbox{$p$ splits completely
in $k/\bbQ$}
\eeq
Thus $p\cO=\fp_1\ldots\fp_r$ where the $\fp_i$
are distinct prime ideals numbered in such a way that $\fp_i$ corresponds
to the embedding $j\circ \tau_i:k\rightarrow \bbQ_p$. Let $\L{i}$
denote the inertia subfield of $\fp_i$ in the extension
$\tK_\infty/k$. We shall need the following consequence
of\refeq{3C}
\begin{lemma}\label{lemma:3A}
For each $i=1,\ldots,r$,
the restriction map $\Gal(\tK_\infty/\L{i})\rightarrow\Gal(\tkinf/k)$
is an isomorphism. In particular $\tL{i}_\infty=\tKinf\ \forall\,i$.
\end{lemma}
\bPf\
The extensions $\L{i}/k$ and  $\tkinf/k$ and  are respectively unramified and
totally ramified at $\fp_i$  so the map is surjective.
Now, $\Gal(\tKinf/\L{i})$ equals $G(\tK_\infty/k)^0_{\fp_i}$.
Completing at any prime above $\fp_i$ it follows from\refeq{3C}
that $\Gal(\tKinf/\L{i})$ is isomorphic to a quotient
of $G(\bbQ_p^{\rm ab}/\bbQ_p)^0$ hence of $\bbZ_p^\times$ by local class field
theory. On the other hand, $\Gal(\tkinf/k)$
is itself isomorphic to $\bbZ_p^\times$ so the surjectivity implies injectivity
and hence both statements of the lemma.\ePf
\begin{lemma}\label{lemma:3B}
Suppose $n\geq v\geq 0$ are integers and $i\in\{1,\ldots,r\}$. Then
\begin{enumerate}
\item\label{part:l3B2} $G(\tL{i}_n/k)^v_{\fp_i}=\Gal(\tL{i}_n/\tL{i}_{v-1})$.
\item\label{part:l3B3} If $v\geq n_2$ then
$G(\tKn/k)_{\fp_i}^v=\Gal(\tKn/\tK_{v-1})$.
\item\label{part:l3B4} $\tL{i}_m=\tK_m$ for all $m\in\bbN$, $m\geq n_2-1$.
\end{enumerate}
\end{lemma}
\bPf\
For part~\ref{part:l3B2}, we complete at a prime of
$\tL{i}_n$ above $\fp_i$. Thus, writing $H$ for the completion of $\L{i}$ at $\fp_i$,
it suffices to show that
$G(\tH_n/\bbQ_p)^v=\Gal(\tH_n/\tH_{v-1})$.
But for every $t\geq -1$, $\tH_t$ is the compositum of
$H$ (which is \emph{unramified} over $\bbQ_p$) with $\bbQ_p(\mu_{p^{t+1}})$
so the proof may be concluded along similar lines to that
of Lemma~\ref{lemma:2B}, \emph{mutatis mutandis}.
Part~\ref{part:l3B3} follows from Lemma~\ref{lemma:2C} on
taking $\fp=\fp_i$ and replacing $K$ by $K(\mu_p)$
so that $K(\mu_p)_t=\tK_t$ for any $t\in\bbN$. (Note that the proof of the
Lemma~\ref{lemma:2B}
does not require $K$ to be totally real. Also, $n_2(K(\mu_p))=n_2(K)=n_2$
since $\ff(K(\mu_p))={\rm l.c.m.}(\ff_0,\ff(k(\mu_p)))={\rm l.c.m.}(\ff_0,p\cO)$.)
For part~\ref{part:l3B4}, we take $v=m+1$ and let $n\rightarrow\infty$ in
parts~\ref{part:l3B2} and~\ref{part:l3B3} to get
$\Gal(\tL{i}_\infty/\tL{i}_{m})=G(\tL{i}_\infty/k)^{m+1}_{\fp_i}$
and $\Gal(\tKinf/\tK_{m})=G(\tKinf/k)_{\fp_i}^{m+1}$. Now apply Lemma~\ref{lemma:3A}.\ePf
\noindent We write $\tGinf$ for $\Gal(\tK_\infty/k)$,
$\tDelta$ for $\Gal(\tK_\infty/\tkinf)$ and $\tGamma{i}$ for
$\Gal(\tKinf/\L{i})$ (see the diagram~(\ref{diag:1})).
Lemma~\ref{lemma:3A} implies the (internal)
direct product decomposition for each $i\in\{1,\ldots,r\}$
\beql{eq:3D}
\tGinf=\tDelta\times\tGamma{i}
\eeq
The appropriate
restriction maps identify $\tDelta$ with the subgroup $\Gal(K/K\cap\tkinf)$
of $G$ and, by\refeq{3D},
with each of the groups $\Gal(\L{i}/k)$ (which are therefore isomorphic
although the $\L{i}$'s may be distinct).
For each $i$ and each $n\geq -1$, we shall write
$\tGamma{i}_n$ for the subgroup $\Gal(\tKinf/\tL{i}_n)$ of $\tGamma{i}$
so that if $n\geq n_2-1$ then $\tGamma{i}_n=\Gal(\tKinf/\tK_n)$ independently of $i$.
Let $\kappa:\tGinf\rightarrow\bbZ_p^\times$ be the
cyclotomic character  so that $\gamma(\zeta)=\zeta^{\kappa(\gamma)}$ for all
$\zeta\in\mu_{p^\infty}:=\bigcup_{n\geq 1}\mu_{p^{n}}$
and $\gamma\in\tGinf$. Then $\kappa$
factors through $\Gal(\tkinf/k)$ and restricts to an
isomorphism $\tGamma{i}\cong\bbZ_p^\times$ for each $i$ taking
$\tGamma{i}_n$ onto $1+\pnpo\bbZ_p$ for each $n\geq 0$.

Now fix $i$ and suppose we are given a sequence
$(v_n)_{n\geq N}$ for some $N\geq 0$ such that,
for all $n\geq N$, $v_n$ lies in $\tL{i}_n$ and satisfies:
\beql{eq:3E}
\mbox{$\ord_\fP(v_n)=0$ for all $\fP|\fp_i$ and
$N_{\tL{i}_{m}/\tL{i}_{n}}v_m=v_{n}\ \forall\,m\geq n$}
\eeq
By taking norms to $\tL{i}_n$ for each $n$ with $N>n\geq 0$
we can, if necessary, extend $\uv$ to a sequence $(v_n)_{n\geq 0}$ (also
denoted $\uv$) with $v_n\in\tL{i}_n$
satisfying\refeq{3E} for all $n\geq 0$. Let $H^{(i)}$ denote
the field $\overline{j\circ\ttau_i(\L{i})}$ (topological closure in $\bbC_p$) which is
a finite, unramified extension of $\bbQ_p$.
For each $\gamma\in\tGinf$ it is clear
that the sequence $j\ttau_i\gamma(\uv):=(j\ttau_i\gamma(v_n))_{n\geq 0}$ lies
in $U_\infty(H^{(i)})$. Suppose $\gamma'\in\tGinf$ fixes $\tL{i}_n$ for some $n\geq 1$.
Since $\tKinf/\tL{i}_n$ is totally ramified at $\fp_i$, we must have
$j\ttau_i\gamma\gamma'=\widehat{\gamma'}j\ttau_i\gamma$
for some $\widehat{\gamma'}\in\Gal(\tH^{(i)}_\infty/\tH^{(i)}_n)$
and the action on $\mu_{p^\infty}$ shows that
$\hat{\kappa}(\widehat{\gamma'})=\kappa(\gamma')=\langle\kappa(\gamma')\rangle$.
Now apply Lemma~\ref{lemma:3Z}
with $H=H^{(i)}$, $\uu=j\ttau_i\gamma\uv$ and
$\sigma=\widehat{\gamma'}$. We deduce easily:
\begin{lemma}\label{lemma:3C}
Fix $t,n\in\bbZ$, $n\geq 1$ and $i$ and $\uv$ as above. As $\gamma$ varies in
$\tGinf$ the value of
$\langle\kappa(\gamma)\rangle^{-t}\check{\delta}_{t,n}(j\ttau_i\gamma\uv)$
depends only on the image of $\gamma$ in $\Gal(\tL{i}_n/k)$.\ePf
\end{lemma}
\noindent
We now define the `higher $p$-adic regulators'. Suppose $\ueps$
is an element of $\cE_\infty(K)$ as defined in the last section.
Then, for each $n\geq N:=\max(n_1,n_2-1)$ we have
$\eps_n\in K_n\subset\tK_n=\tL{i}_{n}$ for $i=1,\ldots,r$
(by Lemma~\ref{lemma:3B}~\ref{part:l3B4}) and for all $n\geq n_2-1$
the norm map $N_{\tL{i}_{m}/\tL{i}_{n}}=N_{\tK_m/\tK_n}$
identifies with $N_{K_m/K_n}$ by restriction
for each $m\geq n$. Therefore $v_n=\eps_n$
satisfies\refeq{3E} for each $n\geq N$ and every $i$.
Clearly, $\Gal(\tKinf/\Kinf)$ fixes $\ueps$ and
lies in the kernel of $\langle\kappa(\,\cdot\,)\rangle$. Therefore,
Lemma~\ref{lemma:3C} implies that for each $i$ and all
$t,n\in\bbZ$, $n\geq n_2(\geq \max(N,1))$,
the quantity $\langle\kappa(\gamma)\rangle^{-t}\check{\delta}_{t,n}(j\ttau_i\gamma\ueps)$
actually depends only on the
image of $\gamma\in\tGinf$ in $\Gal(\tK_n\cap K_\infty/k)=G_n$.
It follows that, for each such $i,t,n$ the following map (depending on $j$)
\displaymapdef{\check{\cL}_{i,t,n}}{\cE_{\infty}(K)}{\cO_{\tH^{(i)}_n}G_n}
{\ueps}
{\sum_{\sigma\in G_n}\langle\kappa(\tilde{\sigma})\rangle^{-t}
\check{\delta}_{t,n}(j\ttau_i\tilde{\sigma}(\ueps))\sigma\inv}
(where $\tilde{\sigma}$ denotes
any lift of $\sigma$ to $\tG_\infty$)
is well-defined and $\bbZ G_\infty$-linear with
$\bbZ\Ginf$ acting on
the group-ring $\cO_{\tH^{(i)}_n}G_n$ through the quotient $\bbZ G_n$.
Suppose $n\geq n_2$. Then for all $\sigma\in G_n$
the above definitions give
\[
j\ttau_i\tilde{\sigma}(\eps_n)^\ast=
j\ttau_i\tilde{\sigma}(\eps_n)^p/j\ttau_i\tilde{\sigma}
(N_{\tL{i}_{n}/\tL{i}_{n-1}}\eps_n)=
j\ttau_i\tilde{\sigma}(\eps_n)^p/j\ttau_i\tilde{\sigma}
(N_{K_n/K_{n-1}}\eps_n)=j\ttau_i\tilde{\sigma}((pe_n)\eps_n)
\]
Putting this together with\refeq{3B.5} and the definition
of $\lambda_{K_n/k,i,p}$, we get
\beql{eq:3F}
\check{\cL}_{i,0,n}(\ueps)=
\frac{1}{p}\lambda_{K_n/k,i,p}((pe_n)\eps_n)=e_n\lambda_{K_n/k,i,p}(\eps_n)
\ \ \ \mbox{for each $i=1,\ldots,r$ and $n\geq n_2$.}
\eeq
Proceeding by analogy with $\lambda_{K_n/k,i,p}$, we extend
each map $\check{\cL}_{i,t,n}$
linearly to $\bbQ\cE_{\infty}(K)$.
Then, for each $n\geq n_2$ and $t\in\bbZ$
we define a unique $\bbQ G_\infty$-linear, (higher) $p$-adic regulator map
$\check{\cR}_{t,n}$ from
${\textstyle\bigwedge_{\bbQ \Ginf}^r}\bbQ\cE_\infty(K)$
to $\bbC_p G_n$
by letting $\check{\cR}_{t,n}(\uu_1\wedge\ldots\wedge\uu_r)=
\det(\check{\cL}_{i,t,n}(\uu_s))_{i,s=1}^r$ for all
$\uu_1,\ldots,\uu_r\in\bbQ\cE_\infty(K)$.
\begin{lemma}\label{lemma:3D}
$\check{\cR}_{0,n}(\ueta)=R_{n,p}(e_n\beta_n(\ueta))$ for all
$\ueta\in{\textstyle\bigwedge_{\bbQ \Ginf}^r}\bbQ\cE_\infty(K)$
and every $n\geq n_2$.
\end{lemma}
\bPf\ By linearity it
suffices to prove the equality when $\ueta$ is of the form
$\bar{\ueps}_1\wedge\ldots\bar{\ueps}_r$ with
$\ueps_i\in\cE_\infty(K)$ for $i=1,\ldots,r$.
But then Equation\refeq{3F} gives
$\check{\cR}_{0,n}(\ueta)=
\det(e_n\lambda_{K_n/k,i,p}(\eps_{s,n}))_{i,s=1}^r=
R_{n,p}(e_n\bar{\eps}_{1,n}\wedge\ldots\wedge e_n\bar{\eps}_{r,n})=
R_{n,p}(e_n\beta_n(\ueta))$.\ePf
\noindent Now let $\ueta$ be
any element of $\bigwedge_{\bbQ \Ginf}^r\bbQ\cE_\infty(K)$
and $m,n$ integers with $n\geq n_2$ and consider the equation
\beql{eq:3G}
\Phi_{n,p}(m)=\frac{2^r}{j(\sqrt{d_k})}\check{\cR}_{1-m,n}(\ueta)
\eeq
The main result of this paper is
\begin{thm}\label{thm:3A}
Suppose that $K$ and $p$ satisfy\refeq{2A},\refeq{2B} and\refeq{3C}.
For a given element $\ueta\in\bigwedge_{\bbQ \Ginf}^r\bbQ\cE_\infty(K)$,
Equation\refeq{3G} holds for all pairs $(m,n)\in\bbZ\times \bbZ_{\geq n_2}$
if and only if it holds for infinitely many such pairs.
\end{thm}
The proof of this theorem will occupy the next section.
Its relevance here comes largely from the following
\begin{cor}\label{cor:3A}
If~\refeq{2A},\refeq{2B} and\refeq{3C} hold and
$\ueta$ demonstrates $P2(K/k,p)$ then it satisfies Equation\refeq{3G}
for all $(m,n)\in\bbZ\times \bbZ_{\geq n_2}$.
\end{cor}
\bPf\ If $\ueta\in\bigwedge_{\bbQ \Ginf}^r\bbQ\cE_\infty(K)$
demonstrates $P2(K/k,p)$ then it satisfies\refeq{3G}
for $m=1$ and all $n\geq n_2$, by
Lemma~\ref{lemma:3D}.\ePf
\section{The Proof~of~Theorem~\ref{thm:3A}}\label{section:pfmainthm}\
To lighten the notation we shall suppress $j$, essentially
regarding $\barbbQ$ simultaneously as a subset of $\bbC$ and of $\bbC_p$,
whenever this cannot lead to confusion.
\begin{figure}
\setlength{\unitlength}{5.5mm}
\beql{diag:1}
\begin{picture}(22,17)(-1,-1)
\put(0.4,0.6){\line(2,3){3.2}}
\put(0,0.6){\line(0,1){3.8}}
\put(0.2,5.5){\line(2,5){1.5}}
\put(0.6,0.4){\line(3,2){22.7}}
\put(3.8,6.4){\line(-1,2){1.4}}
\put(20.4,10.6){\line(2,3){3.2}}
\put(20,10.6){\line(0,1){3.8}}
\put(20.2,15.6){\line(2,5){1.5}}
\put(23.6,16.8){\line(-1,2){1.2}}
\put(0.6,0.3){\line(2,1){18.8}}
\put(2.6,10.3){\line(2,1){5.3}}
\put(0.6,5.3){\line(2,1){6.8}}
\put(4.6,6.3){\line(2,1){18.8}}
\put(11.2,14.6){\line(2,1){9.9}}
\put(8.3,9.75){\line(2,5){1.4}}
\put(9.2,9.6){\line(2,1){10}}
\put(-0.1,0.6){\line(-1,6){1.8}}
\put(-1.1,12.15){\line(6,1){7.9}}
\put(-1.1,12.3){\line(3,1){22.2}}
\put(0.2,0){\makebox(0,0){$k$}}
\put(2,10){\makebox(0,0){$\tK_0$}}
\put(0,5){\makebox(0,0){$K$}}
\put(4.3,5.9){\makebox(0,0){$\tk_0$}}
\put(20.3,10){\makebox(0,0){$k_\infty$}}
\put(22.2,20){\makebox(0,0){$\tK_\infty$}}
\put(20.2,15){\makebox(0,0){$K_\infty$}}
\put(24.1,16){\makebox(0,0){$\tk_\infty$}}
\put(8,9.2){\makebox(0,0){$K_{n_2-1}$}}
\put(10,13.8){\makebox(0,0){$\tL{i}_{n_2-1}=\tK_{n_2-1}$}}
\put(-2,12){\makebox(0,0){$L^{(i)}$}}
\put(10.7,16.5){\makebox(0,0)[br]{$\tGamma{i}$}}
\put(13.3,8.6){\makebox(0,0)[tl]{$\tilde{\Gamma}$}}
\put(10.8,5.1){\makebox(0,0)[tl]{$\Gamma$}}
\put(10.8,9.1){\makebox(0,0)[tl]{$\Gamma$}}
\put(19.85,11.6){\makebox(0,0)[r]{$\Delta$}}
\put(23.2,18.2){\makebox(0,0)[bl]{$\tilde{\Delta}$}}
\end{picture}
\eeq
\end{figure}
In addition to the notations $\tDelta$ and  $\tGamma{i}$ introduced above,
we set $\Delta:=\Gal(K_\infty/k_\infty)$,  $\tiGamma:=\Gal(\tkinf/k)$ and
$\Gamma:=\Gal(k_\infty/k)$ which we identify with the subgroup
$\Gal(\tkinf/\tilde{k}_0)$ of $\tiGamma$ in the obvious way
(see the diagram~(\ref{diag:1})).
To further lighten notation we shall regard the isomorphism
$\tiGamma\cong\bbZ_p^\times$ induced by  $\kappa$ as an identification,
under which $\Gamma$ corresponds to the subgroup $1+p\bbZ_p$.
The restriction to $K_\infty$ gives a surjection
from $\tGinf$ to $\Ginf$ which, with the above identifications,
fits into the exact commuting diagram\refeq{3H}.
\begin{figure}
\setlength{\unitlength}{5.5mm}
\beql{eq:3H}
\begin{picture}(21,8)(0,3)
\put(2,4){\vector(1,0){0.4}}
\put(3.5,4){\vector(1,0){5.5}}
\put(10.7,4){\vector(1,0){4}}
\put(19.2,4){\vector(1,0){0.6}}
\put(2,7){\vector(1,0){0.4}}
\put(3.5,7){\vector(1,0){5.5}}
\put(10.7,7){\vector(1,0){4.8}}
\put(18.5,7){\vector(1,0){1.3}}
\put(3,6.5){\vector(0,-1){2}}
\put(10,6.5){\vector(0,-1){2}}
\put(17,6.5){\vector(0,-1){2}}
\put(17.2,5.5){\makebox(0,0)[l]{$\langle\,\cdot\,\rangle$}}
\put(1.8,4){\makebox(0,0){$1$}}
\put(1.8,7){\makebox(0,0){$1$}}
\put(20,4){\makebox(0,0){$1$}}
\put(20,7){\makebox(0,0){$1$}}
\put(3,4){\makebox(0,0){$\Delta$}}
\put(10,4){\makebox(0,0){$\Ginf$}}
\put(17,4){\makebox(0,0){$\Gamma=1+p\bbZ_p$}}
\put(3,7){\makebox(0,0){$\tDelta$}}
\put(10,7){\makebox(0,0){$\tGinf$}}
\put(17,7){\makebox(0,0){$\tiGamma=\bbZ_p^\times$}}
\end{picture}
\eeq
\end{figure}
The left-hand vertical map here is injective and we regard it henceforth
as an inclusion, \ie\ we identify $\tDelta$
with $\Gal(K_\infty/K_\infty\cap\tkinf)\subset\Delta$.
Lemma~\ref{lemma:3A} shows
that for each $i=1,\ldots,r$ there is a unique injective homomorphism
$\talpha_i:\bbZ_p^\times\rightarrow\tGinf$ with image $\tGamma{i}$ which splits the top
row of\refeq{3H}. Lemma~\ref{lemma:3B}~\ref{part:l3B4}
implies that all the $\talpha_i$ agree on the subgroup
$1+p^{n_2}\bbZ_p$ which they map isomorphically onto $\Gal(\tKinf/\tK_{n_2-1})$.
We fix once and for all an injective homomorphism
$\alpha:1+p\bbZ_p\rightarrow\Ginf$ splitting the bottom
bottom row of\refeq{3H}. For instance, we could choose $\alpha$ to be $\alpha_i$ for
some $i=1,\ldots,r$ where $\alpha_i$ is the composition
of $\talpha_i$ with the inclusion $1+p\bbZ_p\hookrightarrow\bbZ_p^\times$
and the surjection $\tGinf\rightarrow\Ginf$. Thus, whether or not
$\alpha$ actually equals $\alpha_i$ for some $i$, we may clearly assume w.l.o.g. that it \emph{is} chosen to satisfy
\beql{eq:3I}
\alpha(x)=\talpha_i(x)|_{K_\infty}\in\Gal(K_\infty/K_{n_2-1})
\ \ \ \mbox{for all $x\in1+p^{n_2}\bbZ_p$ and all $i\in\{1,\ldots,r\}$}
\eeq
Given such a splitting $\alpha$, the assignment
$\chi\mapsto(\chi|_\Delta,\chi\circ \alpha)$
defines a map
$
\Ginfd\rightarrow\Delta^\dag\times\Gamma^\dag
$
which is clearly bijective. We denote by $\theta\times\psi\in\Ginfd$
the inverse image of $(\theta,\psi)\in\Delta^\dag\times\Gamma^\dag$ under this bijection.

The first step in proving Theorem~\ref{thm:3A} is to break down
Equation\refeq{3G} by characters (linearly extended to the group ring). If
$\ueta\in\bigwedge_{\bbQ \Ginf}^r\bbQ\cE_\infty(K)$ then for any
$(m,n)\in\bbZ\times \bbZ_{\geq n_2}$ we have
\begin{eqnarray}
\lefteqn{\Phi_{n,p}(m)=\frac{2^r}{\sqrt{d_k}}\check{\cR}_{1-m,n}(\ueta)
\Longleftrightarrow}\nonumber\\
&&\chi(\Phi_{n,p}(m))
=
\frac{2^r}{\sqrt{d_k}}\chi(\check{\cR}_{1-m,n}(\ueta))
\ \ \mbox{for all $\chi\in G_\infty^\dag$ factoring through $G_n$}
\label{eq:3J}
\end{eqnarray}
For each character $\chi\in\tGinfd$ we define an integer $n(\chi)\in\bbZ_{\geq -1}$
by
\begin{eqnarray}
n(\chi)&:=&{\rm max}\,\{{\rm ord}_{\fp_i}(\ff(\chi))\,:\,i=1,\ldots,r\}-1\nonumber\\
        &=&{\rm min }\,\{l\geq -1\,:\,
               \mbox{$\chi$ factors through $\Gal(\tL{i}_l/k)$ for $i=1,\ldots,r$}\}\nonumber\\
        &=&{\rm min}\,\{l\geq -1\,:\, (\chi\circ\talpha_i)|_{1+p^{l+1}\bbZ_p}=1\
               \mbox{for $i=1,\ldots,r$}\}\label{eq:3J.5}
\end{eqnarray}
Naturally, `$1+p^0\bbZ_p$' is to be interpreted as
$\bbZ_p^\times$ in the third equality (which then follows from the
definition of $\talpha_i$. The second equality follows from
$\ff(\chi)=\ff(\tKinf^{\ker(\chi)})$, Equation\refeq{2F}
and the fact that $G(\tKinf/k)^{l+1}_{\fp_i}=\Gal(\tKinf/\tL{i}_l)$,
by Lemmas~\ref{lemma:3A} and~\ref{lemma:3B}.)
If $A$ is any quotient group of
$\tGinf$ we shall regard elements of $A^\dag$ as elements of $\tGinfd$ by inflation.
In particular, if $\chi\in\Ginfd$ and
$\psi\in\Gamma^\dag$ then $n(\chi)$ and $n(\psi)$ make sense.
\begin{lemma}\label{lemma:3E}
Suppose $\chi$ is any element of $\Ginfd$ with $\chi=\theta\times\psi$ and
$n\in\bbZ$ with $n\geq n_2-1$.
\begin{enumerate}
\item\label{part:lemma3E1} We have the equivalences:
$(n(\psi)>n)\Leftrightarrow
(\psi(1+p^{n+1}\bbZ_p)\neq\{1\})
\Leftrightarrow (n(\chi)>n)$. In particular $n(\psi)\geq n_2
\Leftrightarrow n(\chi)\geq n_2$,
and in this case $n(\chi)=n(\psi)={\rm min}\{ l\geq n_2\,:\, \psi(1+p^{l+1}\bbZ_p)=\{1\}\}$.
\item\label{part:lemma3E2}
We have the equivalences: $(n(\psi)\leq n)\Leftrightarrow (n(\chi)\leq n)
\Leftrightarrow(\mbox{$\chi$ factors through $G_n$})$.
\end{enumerate}
\end{lemma}
\bPf\ Part~\ref{part:lemma3E1} follows easily from\refeq{3I} and the third equality
in\refeq{3J.5}. The first equivalence
in part~\ref{part:lemma3E2} follows from~\ref{part:lemma3E1}.
For the second, the definition of $n(\chi)$
and Lemma~\ref{lemma:3B}~\ref{part:l3B4} show that
$(n(\chi)\leq n)\Leftrightarrow(\chi$ factors through
$\Gal(\tK_n/k)=\Gal(\tL{i}_n/k)\ \forall\,i)\Leftrightarrow(\chi$ factors
through $\Gal(\tK_n\cap\Kinf/k)=G_n)$.
\ePf
\noindent The following result
establishes the right-hand equation of\refeq{3J} in a
particularly easy case (the proof is deferred).
\begin{prop}\label{prop:3A}
Suppose that $(m,n)\in\bbZ\times \bbZ_{\geq n_2}$ and
$\chi\in\Ginfd$ factors through $G_n$.
If $n(\chi)<n$ then
$\chi(\Phi_{n,p}(m))=0$ and
$\chi(\check{\cR}_{1-m,n}(\ueta))=0$ for all
$\ueta\in\bigwedge_{\bbQ \Ginf}^r\bbQ\cE_\infty(K)$.\ePf
\end{prop}
Let $\Ginfdd=\{\chi\in\Ginfd\,:\,n(\chi)\geq n_2\}$
and $\Gammadd=\{\psi\in\Gamma^\dag\,:\,n(\psi)\geq n_2\}$.
Then Lemma~\ref{lemma:3E}~\ref{part:lemma3E1} implies that the bijection
$
\Ginfd\rightarrow\Delta^\dag\times\Gamma^\dag
$
restricts to a bijection
$
\Ginfdd\rightarrow\Delta^\dag\times\Gammadd
$.
We shall consider functions $\bbZ\times\Gammadd\rightarrow\bbC_p$ and say
that such a function $C$ is \emph{of Iwasawa type} if and only if
there exists a power series $c(T)$ with bounded
coefficients in some finite extension of $\bbQ_p$ such that
$C(m,\psi)=c((1+p)^m\psi(1+p)\inv-1)$ for every $(m,\psi)\in\bbZ\times\Gammadd$.
(The series converges because $\psi$ takes values in $\mu_{p^\infty}$.)
We write $\IwaZG$ for the set of all such functions, considered as a $\bar{\bbQ}_p$-algebra
under pointwise operations.
By the Uniqueness Principle, $c(T)$ is unique given $C$, if it exists. Note also
that the r\^{o}le of $(1+p)$ in this definition could be played by
any topological generator of $1+p\bbZ_p$: the power series $c$ would
change but the set $\IwaZG$ would not.
For each $\theta\in\Delta^\dag$ we define
a function
\displaymapdef{C_{\theta}}{\bbZ\times\Gammadd}{\bbC_p}{(m,\psi)}
{g(\hchi)\inv\langle N\ff'(\chi)d_k\rangle^{1-m}\chi\left(\Phi_{n(\psi),p}(m)\right)
\ \mbox{where $\chi=\theta\times\psi\in\Ginfdd$}}
Note that here and in what follows, the set
$T$ is implicit and equal to $T_p$. Therefore the notations
$\ff'(\chi)$, $\ff'(K)$, $\ff'(K_n)$ \etc\
represent the \emph{prime-to-$p$} parts of $\ff(\chi)$, $\ff(K)$,
$\ff(K_n)$ \etc\
Similarly $\cD'$ denotes the prime-to-$p$ part of
$\cD$, which equals $\cD$ by
Condition\refeq{3C} (or\refeq{2C}). In particular $N\ff'(\chi)d_k=N(\ff'(\chi)\cD')$.
For each $\theta\in\Delta^\dag$ and $\ueta\in\WedQGEinf$ we define
\displaymapdef{C_{\theta,\ueta}}{\bbZ\times\Gammadd}{\bbC_p}{(m,\psi)}
{\displaystyle\frac{2^r}{\sqrt{d_k}}
g(\hchi)\inv\langle N\ff'(\chi)d_k\rangle^{1-m}
\chi\left(\check{\cR}_{1-m,n(\psi)}(\ueta)\right)
\ \mbox{where $\chi=\theta\times\psi\in\Ginfdd$}}
\begin{prop}\label{prop:3B}
For each $\theta\in\Delta^\dag$ the function $C_\theta$ lies in $\IwaZG$
\ePf
\end{prop}
\begin{prop}\label{prop:3C}
For each $\theta\in\Delta^\dag$ and every $\ueta\in\WedQGEinf$
the function $C_{\theta,\ueta}$ lies in $\IwaZG$
\ePf
\end{prop}
The proofs of Propositions~\ref{prop:3B} and~\ref{prop:3C} are also deferred
while we show how, along with Proposition~\ref{prop:3A}, they imply
Theorem~\ref{thm:3A}.
First note that for any choice of $(m,\psi)$ in
$\ZGdd$, $\theta$ in $\Deltad$ and $\ueta\in\WedQGEinf$,
the character $\chi=\theta\times\psi$ factors through
$G_{n(\chi)}=G_{n(\psi)}$ (by Lemma~\ref{lemma:3E}) and
we have the equivalences:
\begin{eqnarray}
&&
\chi\left(\Phi_{n(\psi),p}(m)\right)=\frac{2^r}{\sqrt{d_k}}
\chi\left(\check{\cR}_{1-m,n(\psi)}(\ueta)\right)\ \
\mbox{(with $\chi=\theta\times\psi$)}
\label{eq:3K}\\
&\Leftrightarrow&C_{\theta}(m,\psi)=C_{\theta,\ueta}(m,\psi)
\nonumber\\
&\Leftrightarrow&c_{\theta}((1+p)^m\psi(1+p)\inv-1)=
c_{\theta,\ueta}((1+p)^m\psi(1+p)\inv-1)
\label{eq:3M}
\end{eqnarray}
where $c_{\theta}$ and $c_{\theta,\ueta}$
are the power series whose existence
is implied by the statements of
Propositions~\ref{prop:3B} and~\ref{prop:3C} respectively.
Suppose that $\ueta_0$ is an element of $\WedQGEinf$ such that\refeq{3G}
is satisfied
for an infinite subset $Y$, say, of $\ZZnt$ and let $Y'$ be the
inverse image of this set under the map
$\ZGdd\rightarrow\ZZnt$ sending $(m,\psi)$ to $(m,n(\psi))$.
Lemma~\ref{lemma:3E} shows that this map is surjective so, in
particular, $Y'$ is also infinite. But the equivalence\refeq{3J}
shows that Equation\refeq{3K} holds for each $(m,\psi)\in Y'$, with
$\ueta=\ueta_0$ and any $\theta\in\Deltad$. Hence the same is true of Equation\refeq{3M}.
Since the map $(m,\psi)\mapsto((1+p)^m\psi(1+p)\inv-1)$ is injective on $Y'$,
it follows that for any $\theta\in\Deltad$, the power series
$c_\theta$ and $c_{\theta,\ueta_0}$ agree at infinitely
many $a\in\bbC_p$ with $|a|_p<1$. The Uniqueness Principle
therefore implies that $c_\theta=c_{\theta,\ueta_0}$ for all
$\theta\in\Deltad$. Hence Equation\refeq{3M} holds with
$\ueta=\ueta_0$ for any $\theta\in\Deltad$ and $(m,\psi)\in \ZGdd$.
Hence the same is true of Equation\refeq{3K}.
In other words, for any $n\geq n_2$, the condition
\beql{eq:3N}
\chi\left(\Phi_{n,p}(m)\right)=\frac{2^r}{\sqrt{d_k}}
\chi\left(\check{\cR}_{1-m,n}(\ueta_0)\right)\ \ \ \mbox{$\forall\,m\in\bbZ$}
\eeq
is satisfied by all $\chi\in\Ginfd$ with $n(\chi)=n$.
But by Lemma~\ref{lemma:3E} again, any other character
$\chi\in\Ginfd$ factoring through $G_n$ must have
$n(\chi)<n$, in which case Proposition~\ref{prop:3A}
implies that it satisfies\refeq{3N} trivially.
Thus the equivalence\refeq{3J} shows that Equation\refeq{3G} must hold
with $\ueta=\ueta_0$ for all $m\in\bbZ$ and all $n\geq n_2$, as required.

In order to prove Proposition~\ref{prop:3B}, we
first evaluate the term
$\chi\left(\Phi_{n(\psi),p}(m)\right)$ in the definition of
$C_\theta(m,\psi)$.
If $L$ is an infinite abelian extension of $k$ and $\fq$ a prime ideal
of $\cO$ then for each $v\in[-1,\infty)$ we denote by $G(L/k)^v_\fq$
the $v$th ramification subgroup (in the upper numbering)
of $\Gal(L/k)$, namely
the inverse limit of the $G(L'/k)^v_\fq$ as $L'/k$ runs through all finite
subextensions of $L/k$. If $\fq\ndiv p$ and $v\geq 0$ then the fact that
$\kinf/k$ is unramified at $\fq$ implies
that $G(\Kinf/k)^v_\fq$ is contained in
$\Delta$ and so maps isomorphically onto $G(K/k)^v_\fq$.
This implies in particular that for any
$\theta\in\Deltad$ it makes sense to define
an ideal $\ff'(\theta)$ of $\cO$ \emph{which is prime to $p$} by
\beql{eq:3N.5}
\ord_\fq(\ff'(\theta))=
\min\{v\in\bbN:G(\Kinf/k)^v_\fq\subset\ker (\theta)\}
\ \ \ \mbox{for all primes $\fq$
of $\cO$ not dividing $p$.}
\eeq
Let $\chi=\theta\times\psi\in\Ginfd$ so in particular
$\theta=\chi|_\Delta$. Since $\ff(\chi)=\ff(\Kinf^{\ker(\chi)})$ and
$G(\Kinf/k)^v_\fq\subset\Delta$ for each $\fq\ndiv p$ and all $v\geq 0$,
it follows from\refeq{2F} and\refeq{3N.5}
that
\beql{eq:3N.6}
\ff'(\chi)=\ff'(\theta)
\eeq
(where $\ff'(\chi)$ denotes the prime-to-$p$ part of $\ff(\chi)$).
If $\chi$ lies in $\Ginfdd$ then
$n(\psi)\geq n_2$ and $\chi$ factors through $G_{n(\psi)}$.
The definition of $n(\chi)$
shows that $\ord_{\fp_i}(\ff(\chi))=n(\chi)+1=n(\psi)+1$ for some $i$,
hence for all $i=1,\ldots,r$, by\refeq{2H} with $n=n(\psi)$.
Therefore
\beql{eq:3N.75}
\chi=\theta\times\psi\in\Ginfdd\Longrightarrow
\ff(\chi)=p^{n(\psi)+1}\ff'(\theta)
\eeq
On the other hand, Equations\refeq{2G} and\refeq{xxtra} show
that $\ff(K_{n(\psi)})=
\ff_{n(\psi)}=\ff'(K)p^{n(\psi)+1}\cO$. Thus,
applying Proposition~\ref{prop:1A}
(replacing $K$ by $K_{n(\psi)}$ and $T$ by $T_p$)
we find that $\fh_0=\cO$ and (with a slight change of notation)
Equation\refeq{1H} gives
\begin{eqnarray}
\lefteqn{
\theta\in\Deltad,
(m,\psi)\in\ZGdd,\ \chi=\theta\times \psi\Longrightarrow}
\nonumber\\
&&\chi\left(\Phi_{n(\psi),p}(m)\right)=g(\hchi)\langle
N\ff'(\chi) d_k\rangle^{m-1}
\prod_{\fq\ \mathit{prime}\atop \fq|\ff'(K),\ \fq\ndiv\ff'(\theta)}
B_{\theta,\fq}(m,\psi)L_p(m,\hchi)
\label{eq:3O}
\end{eqnarray}
where, for
each $\theta\in\Deltad$ and prime ideal $\fq\ndiv p\ff'(\theta)$
of $\cO$ we define the function
\displaymapdef{B_{\theta,\fq}}{\bbZ\times\Gammadd}{\bbC_p}{(m,\psi)}
{\displaystyle 1-\frac{\hchi\inv(\fq)}{\langle N\fq\rangle ^{1-m}}
\ \ \mbox{where $\chi=\theta\times\psi\in\Ginfdd$}}
For any $x\in\bbZ_p^\times$ we set $\ell(x):=\log_p(x)/\log_p(1+p)\in \bbZ_p$
so that $(1+p)^{\ell(x)}=\langle x \rangle$. We shall need two simple lemmas.
\begin{lemma}\label{lemma:3F}
For each $\theta\in\Delta^\dag$ and $\fq\ndiv p\ff'(\theta)$
the function $B_{\theta,\fq}$ lies in $\IwaZG$
\end{lemma}
\bPf\ Let $L_\fq$ be the inertia subfield of $\fq$ in $\Kinf/k$
so that $L_\fq\supset \kinf$ and the decomposition group
$G(L_\fq/k)_\fq^{-1}$ is (infinite) procyclic with canonical topological
generator $\sigma_\fq$, say, given by the Frobenius at
primes above $\fq$ in every finite sub-extension. We fix a lift
$\tilde{\sigma}_\fq\in \Ginf$ of $\sigma_\fq$ and let
$\delta_\fq=\tilde{\sigma}_\fq\alpha(\tilde{\sigma}_\fq|_{\kinf})\inv\in\Delta$.
Note that $\tilde{\sigma}_\fq|_{\kinf}\in\Gal(\kinf/k)$ identifies
with $\langle N\fq \rangle\in 1+p\bbZ_p$ under $\kappa$.
Let $\theta\in\Deltad$,
$(m,\psi)\in\ZGdd$ and $\chi=\theta\times\psi$.
The condition $\fq\ndiv p\ff'(\theta)=p\ff'(\chi)$ means that
$\chi$ factors through $\Gal(L_\fq/k)$ and $\hchi(\fq)=\chi(\sigma_\fq)=
\chi(\tilde{\sigma}_\fq)=\theta(\delta_\fq)\psi(\tilde{\sigma}_\fq|_{\kinf})=
\theta(\delta_\fq)\psi(\langle N\fq\rangle)$. Thus
$B_{\theta,\fq}(m,\psi)$ can be written as
$1-\theta(\delta_\fq)\inv\langle N\fq \rangle\inv
\left((1+p)^m\psi(1+p)\inv\right)^{\ell(N\fq)}$. The lemma now
follows from the fact that $(1+T)^{\ell(N\fq)}\in\bbZ_p[[T]]$.\ePf
\noindent Suppose that $U$ is an open and closed subset of $\bbZ_p^\times$
and $\mu\in{\rm Meas}(U,L)$ for some finite extension $L$ of $\bbQ_p$.
We shall write $I_{U,\mu}$ for the function
$\ZGdd\rightarrow\bbC_p$ sending $(m,\psi)$ to
$\int_{x\in U}\langle x\rangle^{-m}\psi(\langle x\rangle)d\mu$.
\begin{lemma}\label{lemma:3G}
For $U$, $\mu$ as above, the function
$I_{U,\mu}$ lies in $\IwaZG$.
\end{lemma}
\bPf\ If $x\in\bbZ_p^\times$ then the Binomial Theorem gives
\[
\langle x \rangle^{-m}\psi(\langle x\rangle)=((1+p)^{m}\psi(1+p)\inv)^{-\ell(x)}=
\sum_{t=0}^{\infty}
\left(
\begin{array}{c}
-\ell(x)\\
t
\end{array}
\right)
((1+p)^{m}\psi(1+p)\inv-1)^t
\]
where the series converges uniformly as a function of $x\in\bbZ_p^\times$
with values in $\bbZ_p[\psi]$.
Therefore, $I_{U,\mu}(m,\psi)$ may be written as $\sum_{t=0}^{\infty}
c_{U,\mu,t}((1+p)^{m}\psi(1+p)\inv-1)^t$ where $c_{U,\mu,t}=\int_{x\in U}
\left(
\begin{array}{c}
-\ell(x)\\
t
\end{array}
\right)
d\mu\in L$ is bounded as a function of $t$.\ePf
\noindent For each $\theta\in\Deltad$, Equation\refeq{3O} implies that for
every $(m,\psi)\in\ZGdd$
\beql{eq:3P}
C_\theta(m,\psi)=\prod_{\fq\ \mathit{prime}\atop \fq|\ff'(K),\ \fq\ndiv\ff'(\theta)}
B_{\theta,\fq}(m,\psi)L_p(m,\widehat{\theta\times\psi})
\eeq
The proof of Proposition~\ref{prop:3B} now splits into two cases according as
$\theta$ is or is not the trivial character.\\
\underline{\bf Case 1: $\theta\neq 1$} In this case the Proposition
will follow from\refeq{3P} and Lemma~\ref{lemma:3F} once we have proven
\beql{eq:3P.5}
\mbox{$\theta\neq 1\Longrightarrow$ the map
$\ZGdd\ni (m,\psi)\mapsto L_p(m,\widehat{\theta\times\psi})$ is of Iwasawa type.}
\eeq
We now deduce this implication from the results of~\cite[\S 4]{Rib}.
To do this properly, we first need to
connect our notations and set-up with those of~\emph{ibid.}
For Ribet's ideal `$f$' we take $\ff'(\theta)p\cO$ so that
Ribet's group `$G$' identifies with $\cG_{\theta,\infty}
:=\Gal(\cK_{\theta,\infty}/k)$ where
$\cK_{\theta,\infty}$ is the union over $n$
of the \emph{strict} ray-class
fields over $k$ modulo $\ff'(\theta)p^{n+1}\cO$.
For every $\psi\in\Gammad$ the element of $\Ginfd$ which we denote $\theta\times \psi$
factors through $\Gal(K_\infty^{\ker(\theta)}/k)$ and since
$\cK_{\theta,\infty}$ clearly contains $K_\infty^{\ker(\theta)}$, we may regard
$\theta\times \psi$ as an element of $\cG^\dag_{\theta,\infty}$.
In particular, we shall take the character of $\cG_{\theta,\infty}$
denoted `$\varepsilon$' by Ribet to be $\theta\times 1$.
Our group $\Gamma$ (identified with $1+p\bbZ_p$)
coincides with Ribet's. His splitting of
$\cG_{\theta,\infty}$ as `$\Gamma\times A$' (with $A=\Gal(\cK_{\theta,\infty}/\kinf)$)
depends on the choice of a splitting homomorphism
$\tilde{\alpha}$, say, from $\Gamma$ to $\cG_{\theta,\infty}$.
We may (and shall) choose $\tilde{\alpha}$ to be a lift
of the injective homomorphism $\Gamma\rightarrow\Gal(K_\infty^{\ker(\theta)}/k)$
obtained by composing our map $\alpha:\Gamma\rightarrow\Ginf$ with
the natural surjection $\Ginf\rightarrow\Gal(K_\infty^{\ker(\theta)}/k)$.
Any $\psi\in\Gamma^\dag$ inflates to
$\tilde{\psi}\in\cG^\dag_{\theta,\infty}$ which is trivial on $A$
and the above choice means that element of $\cG^\dag_{\theta,\infty}$
denoted `$\tilde{\psi}\varepsilon$' in Ribet's notation coincides
with $\theta\times\psi$.
Moreover, if $\psi$ lies in $\Gammadd$ then\refeq{3N.75} shows that
the prime factors of $\ff(\theta\times\psi)$ are those
of $\ff'(\theta)p\cO$, so that the complex function which would be denoted
`$L(s,\tilde{\psi}\varepsilon)$' by Ribet in this situation
(See pp. 179 \emph{and}~186 of~\cite{Rib}) is exactly equal to our
primitive $L$-function
$L(s,\widehat{\theta\times\psi})$ (see~\eg~\cite[\S2]{zetap1}).
Interpolating this equality $p$-adically for $s$ in the dense
subset $\cM(p)$ of $\bbZ_p$, it follows that (in an obvious notation)
\beql{eq:3Q}
L_{p,{\rm Ribet}}(m,\tilde{\psi}\varepsilon)=L_p(m,\widehat{\theta\times\psi})
\ \ \ \forall\,(m,\psi)\in\ZGdd
\eeq
Now, $\theta\neq 1$ implies that
`$\varepsilon$ is non-trivial on $A$'
and so, taking Ribet's generator `$\gamma$' of $\Gamma$ to be $1+p$, his
statements~(4.9) and (4.10) imply in this case that there exists a power series
`$F_\varepsilon(T)$' with coefficients in $R=\bbZ_p[\theta]$
such that $L_{p,{\rm Ribet}}(s,\tilde{\psi}\varepsilon)=F_\varepsilon(\psi(1+p)(1+p)^{1-s}-1)$
for all $s\in\bbZ_p$ and $\psi\in\Gammadd$. Hence, if $\theta\neq 1$, Equation\refeq{3Q}
gives $L_p(m,\widehat{\theta\times\psi})=
c_\theta((1+p)^m\psi(1+p)\inv-1)$ for all $(m,\psi)\in\ZGdd$,
where $c_\theta(T)=F_\varepsilon((1+p)(1+T)\inv -1)\in\bbZ_p[\theta][[T]]$,
thus proving\refeq{3P.5}.\\
\underline{\bf Case 2: $\theta=1$}. The comparison of our set-up
with Ribet's, and in particular Equation\refeq{3Q}, goes
through exactly as in Case~1. Since
Ribet's `$\varepsilon$' is now the trivial character of $\cG_{\theta,\infty}$ however,
his statements~(4.9) and (4.10)
do not apply directly. On the other hand, since now $\ff'(\theta)=\cO$,
the product in\refeq{3P} extends over the set of all prime divisors $\fq$
of $\ff'(K)$ and
Condition\refeq{2B} says that this set is \emph{non-empty}. So if we fix
$\fq_0$ dividing $\ff'(K)$ then
Proposition~\ref{prop:3B} will follow
from Lemma~\ref{lemma:3F} and Equations\refeq{3P} and\refeq{3Q} provided
we can show that
\beql{eq:3R}
\mbox{the map
$\ZGdd\ni (m,\psi)\mapsto B_{1,\fq_0}(m,\psi)L_{p,{\rm Ribet}}(m,\tilde{\psi})$
is of Iwasawa type.}
\eeq
The field $\cK_{\theta,\infty}=\cK_{1,\infty}$ is now the union over $n$
of the strict ray-class
fields over $k$ modulo $p^{n+1}\cO$.
We (temporarily) write $\pi$ for the natural surjection from
$\cG_{1,\infty}=\Gal(\cK_{1,\infty}/k)$ to $\Gamma$
(identified with $1+p\bbZ_p$)
and $\sigma_0$ for the element of $\cG_{1,\infty}$ which is the
limit of the Frobenius at primes above $\fq_0$ in finite subextensions.
For each $(m,\psi)\in\ZGdd$ we
substitute  $m$ for $s$,
$\sigma_0$ for `$c$' and $\tilde{\psi}$ for `$\epsilon$' in Ribet's~(4.6) to obtain
(in our notation)
\[
(1-\tilde{\psi}(\sigma_0)\langle N\fq_0\rangle^{1-m})L_{p,{\rm Ribet}}(m,\tilde{\psi})=
\int_{y\in\cG_{1,\infty}}\pi(y)^{1-m}\tilde{\psi}(y)d\lambda_{\sigma_0}
\]
where $\lambda_{\sigma_0}$ is a (bounded) measure on
the group $\cG_{1,\infty}$ with
values in $\bbQ_p$ as constructed by Ribet. We use $\pi$ to push forward
the measure $\pi(y)\lambda_{\sigma_0}(y)$
on $\cG_{1,\infty}$ to one on $1+p\bbZ_p$ giving an element
$\mu_{\sigma_0}$, say, of ${\rm Meas}(1+p\bbZ_p,\bbQ_p)$. Since also
$\tilde{\psi}(\sigma_0)=\widehat{(1\times\psi)}(\fq_0)$ by definition,
the last equation may be rewritten
\[
\frac{B_{1,\fq_0}(m,\psi)}{B_{1,\fq_0}(m,\psi)-1}
L_{p,{\rm Ribet}}(m,\tilde{\psi})=
\int_{x\in \Gamma} x^{-m}\psi(x)d\mu_{\fq_0}=
I_{1+p\bbZ_p,\mu_{\fq_0}}(m,\psi)
\]
Multiplying both sides by $B_{1,\fq_0}(m,\psi)-1$ and using
Lemmas~\ref{lemma:3F}
and~\ref{lemma:3G} establishes\refeq{3R}.
This completes the proof of Proposition~\ref{prop:3B}.

We now introduce some notation and two lemmas that will be used
in the proofs of Propositions~\ref{prop:3A} and~\ref{prop:3C}.
For any $n\in\bbN$ and any ($p$-adic) character
$\tpsi\in(\bbZ_p^\times)^\dag$ with $1+p^{n+1}\bbZ_p\subset\ker{\tpsi}$
we choose a set $R_n$
of representatives for $\bbZ_p^\times$ modulo
$1+p^{n+1}\bbZ_p$ and
define the Gauss sum
\[
g_n(\tpsi):=\sum_{u\in R_n}
\zeta_n^u\tpsi(u)\inv\in\bbC_p
\]
which is clearly independent of the choice of $R_n$.
\begin{lemma}\label{lemma:3H}
Suppose that $n\geq 1$ and $1+p^n\bbZ_p\subset\ker{\tpsi}$. Then
$g_n(\tpsi)=0$.
\end{lemma}
\bPf\ If $R_n$ is a set of representatives as above then,
since $n\geq 1$, so is $R'_n:=\{u+p^n:u\in R_n\}$.
Since $\tpsi(u)=\tpsi(u+p^n)$ we get
\[
\zeta_0g_n(\tpsi)=
\zeta_n^{p^n}g_n(\tpsi)=
\sum_{u\in R_n}
\zeta_n^{u+p^n}\tpsi(u)\inv=
\sum_{u'\in R'_n}
\zeta_n^{u'}\tpsi(u')\inv=
g_n(\tpsi)
\]
but $\zeta_0\neq 1$, so $g_n(\tpsi)=0$.\ePf
\noindent
For $i=1,\ldots,r$ we let
$\trho_i$ be the unique element of
$\tDelta$ lifting the Frobenius element above $\fp_i$ in
$\Gal(L^{(i)}/k)$ (see\refeq{3D}).
\begin{lemma}\label{lemma:3I}
Suppose $i\in\{1,\ldots,r\}$, $\theta\in\Deltad$ and
$\ueps\in\cE_\infty(K)$. Then there exists a measure
$\mu_{i,\theta,\ueps}\in {\rm Meas}(\bbZ_p^\times,H^{(i)}(\theta))$
depending only on $i$, $\theta$ and $\ueps$ and
such that for any $n\in\bbZ_{\geq n_2}$ and for all
$(m,\psi)\in\bbZ\times\Gammad$ satisfying  $n(\psi)\leq n$, we have
\beql{eq:3S}
\chi\left( \check{\cL}_{i,1-m,n}(\ueps)\right)=
g_n(\chi\circ\talpha_i)\theta(\trho_i)^{-n}
\int_{x\in\bbZ_p^\times}\langle x\rangle^{1-m}(\chi\circ\talpha_i)(x)d\mu_{i,\theta,\ueps}
\eeq
where $\chi=\theta\times\psi$.
\end{lemma}
\rem\ Lemma~\ref{lemma:3E}~\ref{part:lemma3E2}) implies that
$\chi$ factors through $G_n$ so the L.H.S. of\refeq{3S} makes sense. It
also implies that $n(\chi)\leq n$, so $g_n(\chi\circ\talpha_i)$ makes sense, by\refeq{3J.5}.
The proof will show that $\mu_{i,\theta,\ueps}$ actually
depends only on $\theta|_\tDelta$.\vspace{1ex}\\
\bPf\
As explained above, the conditions of the Lemma allow us
to regard $\chi$ as a character of $\Ginf$ (hence of $\tGinf$)
factoring through $G_n$. The definition of $\check{\cL}_{i,1-m,n}$
therefore gives $\chi\left( \check{\cL}_{i,1-m,n}(\ueps)\right)=
\sum\langle\kappa(\tilde{\sigma})\rangle^{m-1}
\check{\delta}_{1-m,n}(\ttau_i\tilde{\sigma}(\ueps))
\chi(\tilde{\sigma})\inv$ where $\tilde{\sigma}$ runs through any set of representatives
for $\tGinf$ modulo $\Gal(\tKinf/K_n)$.
Choose a set $R_n$ of representatives for $\bbZ_p^\times$ modulo
$1+p^{n+1}\bbZ_p$. Then it follows from\refeq{3D} and the definition
of $\talpha_i$ that
$\{\talpha_i(u)\td\,:\,u\in R_n,\ \td\in\tDelta\}$ is a set of representatives of
$\tGinf$ modulo $\talpha_i(1+p^{n+1}\bbZ_p)=\tilde{\Gamma}^{(i)}_n$ which equals
$\Gal(\tKinf/\tK_n)$ since $n\geq n_2-1$. It follows that
\begin{eqnarray*}
\chi\left( \check{\cL}_{i,1-m,n}(\ueps)\right)&=&
[\tK_n:K_n]\inv\sum_{u\in R_n}\sum_{\td\in\tDelta}\langle u\rangle^{m-1}
\check{\delta}_{1-m,n}(\ttau_i\talpha_i(u)\td\ueps)
(\chi\circ\talpha_i)(u)\inv\theta(\td)\inv\\
&=&a\inv\sum_{u\in R_n}\sum_{\td\in\tDelta}\langle u\rangle^{m-1}
\phi^{-n}\check{D}^{1-m}h^\ast(\ttau_i\talpha_i(u)\td\ueps;\pi_n)
(\chi\circ\talpha_i)(u)\inv\theta(\td)\inv\\
\end{eqnarray*}
where $a:=[\tK_0:K]=[\tK_n:K_n]\ \forall n$. The action on $\mu_{p^\infty}$
shows that $\ttau_i\talpha_i(u)=\hat{\kappa}\inv(u)\ttau_i$.
(Recall that $\hat{\kappa}$ denotes the Lubin-Tate isomorphism
from $\Gal(\tilde{H}_\infty^{(i)}/H^{(i)})$ to $\bbZ_p^\times$.)
Therefore Lemma~\ref{lemma:3Z} gives
$\check{D}^{1-m}h^\ast(\ttau_i\talpha_i(u)\td\ueps;\pi_n)=
\langle u\rangle^{1-m}\check{D}^{1-m}h^\ast(\ttau_i\td\ueps;\zeta_n^u-1)$ and so
\[
\chi\left( \check{\cL}_{i,1-m,n}(\ueps)\right)=
\sum_{u\in R_n}(\chi\circ\talpha_i)(u)\inv
\check{D}^{1-m}
\left.\left[a\inv\sum_{\td\in\tDelta}\theta(\td)\inv\phi^{-n}
h^\ast(\ttau_i\td\ueps;X)\right]\right|_{X=\zeta_n^u-1}
\]
Now for any $\underline{v}\in U_\infty(H^{(i)})$ it is clear that
$\phi^{-n}g(\underline{v};X)=g(\phi^{-n}\underline{v};X)$ and hence that the
same equality holds with $g^\ast$ in place of $g$ on both sides,
hence the same with $h^\ast$ in place of $g$.
Moreover, it is easy to see that $\phi^{-n}\ttau_i=\ttau_i\tilde{\rho}_i^{-n}$.
Therefore
\[
a\inv\sum_{\td\in\tDelta}\theta(\td)\inv\phi^{-n}
h^\ast(\ttau_i\td\ueps;X)=
a\inv\sum_{\td\in\tDelta}\theta(\td)\inv
h^\ast(\ttau_i\tilde{\rho}_i^{-n}\td\ueps;X)=
\theta(\tilde{\rho}_i)^{-n}F(i,\theta,\ueps;X)
\]
where $F(i,\theta,\ueps;X)$ is defined to be the power series
$a\inv\sum_{\td\in\tDelta}\theta(\td)\inv h^\ast(\ttau_i\td\ueps;X)$,
which lies in $\cA^\ast(H^{(i)}(\theta))$.
(Since $p$ is odd, $a$ divides $p-1$,
from which it follows easily that $F(i,\theta,\ueps;X)$ actually has
coefficients in the ring of valuation integers of $H^{(i)}(\theta|_\tDelta)$.)
Thus $F(i,\theta,\ueps;X)$ corresponds
to an element of ${\rm Meas}(\bbZ_p^\times,H^{(i)}(\theta))$
which we denote $\mu_{i,\theta,\ueps}$. As explained above,
the formalism of the power series/measure
correspondence implies that
$\check{D}^{1-m}F(i,\theta,\ueps;X)$ corresponds
to the measure $\langle x\rangle^{1-m}\mu_{i,\theta,\ueps}(x)$
and also that the value of a power series
at $\zeta-1$ (for any $\zeta\in\mu_{p^\infty}$) is equal to the integral of the
function $x\mapsto\zeta^x$ against the corresponding measure.
Putting all this together,
we obtain
\begin{eqnarray*}
\chi\left(\check{\cL}_{i,1-m,n}(\ueps)\right)&=&
\theta(\tilde{\rho}_i)^{-n}
\sum_{u\in R_n}(\chi\circ\talpha_i)(u)\inv
\int_{x\in\bbZ_p^\times}\zeta_n^{ux} \langle x\rangle^{1-m}d\mu_{i,\theta,\ueps}\\
&=&
\theta(\tilde{\rho}_i)^{-n}
\int_{x\in\bbZ_p^\times}
\left(\sum_{u\in R_n}\zeta_n^{ux}(\chi\circ\talpha_i)(ux)\inv\right)
\langle x\rangle^{1-m}(\chi\circ\talpha_i)(x)d\mu_{i,\theta,\ueps}
\end{eqnarray*}
But, for any $x\in\bbZ_p^\times$, the set $\{xu:u\in R_n\}$ is clearly
another set of representatives of $\bbZ_p^\times$ modulo $1+p^{n+1}\bbZ_p$
and Equation\refeq{3S} follows.\ePf
\noindent We can now prove Proposition~\ref{prop:3A}.
The general form of an element $\ueta\in\WedQGEinf$ is as in\refeq{2E.25}.
The second equality of Proposition~\ref{prop:3A}
will follow by linearity of $\check{\cR}_{1-m,n}$ if we can
prove it in the special case
\beql{eq:3S.5}
\ueta=\ueps_1\wedge\ldots\wedge\ueps_r
\ \ \ \mbox{with $\ueps_1,\ldots,\ueps_r\in\cE_\infty(K)$}
\eeq
so that, for any
$m\in\bbZ$, $n\in\bbZ_{\geq n_2}$ and $\chi\in\Ginfd$ factoring through $G_n$ we may write
\beql{eq:3T}
\chi(\check{\cR}_{1-m,n}(\ueta))=\det(\chi(\check{\cL}_{i,1-m,n}(\ueps_l)))_{i,l=1}^r
\eeq
Suppose $n(\chi)<n$.
Then on the one hand $\chi=\theta\times\psi$ with $n(\psi)<n$
(by Lemma~\ref{lemma:3E}~\ref{part:lemma3E2}) so each term
$\chi(\check{\cL}_{i,1-m,n}(\ueps_l))$ in\refeq{3T}
may be calculated by means of\refeq{3S}. On the other hand
Equation\refeq{3J.5} implies that for all $i$ we have
$1+p^n\bbZ_p\subset \ker(\chi\circ\talpha_i)$, hence $g_n(\chi\circ\talpha_i)$ vanishes
by Lemma~\ref{lemma:3H}. Thus $\chi(\check{\cL}_{i,1-m,n}(\ueps_l))=0$ for all
$i,l,m$, hence the equality $\chi(\check{\cR}_{1-m,n}(\ueta))=0$.
For the other equality in Proposition~\ref{prop:3A},
we note that since $n(\chi)<n$,
Equations\refeq{3J.5} and\refeq{2G} imply that for some $i$ one has
$\ord_{\fp_i}(\ff(\chi))<n+1=\ord_{\fp_i}(\ff_n)$, \ie\
$\ff(\chi)|\fp_i\inv\ff_n$.
(In fact, this must hold for for all $i$ by\refeq{2H}, since $n\geq n_2$).
The equality $\chi(\Phi_{n,p}(1))=0$ therefore follows
from\refeq{1J} (\cf\ the proof of Proposition~\ref{prop:2A}). This
completes the proof of Proposition~\ref{prop:3A}.

For Proposition~\ref{prop:3C} we introduce two new functions as follows.
For each $\theta\in\Deltad$, $i=1,\ldots,r$ and
$\ueps\in\cE_\infty(K)$ we define
$E_\theta,\ D_{i,\theta,\ueps}:\bbZ\times\Gammadd\rightarrow\bbC_p$
by setting $E_\theta(m,\psi)=\langle N\ff'(\theta)d_k\rangle^{1-m}
\psi(\langle N\ff'(\theta)d_k\rangle)$ and
$D_{i,\theta,\ueps}(m,\psi)=
\int_{x\in\bbZ_p^\times}\langle x\rangle^{1-m}
(\chi\circ\talpha_i)(x)d\mu_{i,\theta,\ueps}$ with
$\chi=\theta\times\psi\in\Ginfdd$ and $\mu_{i,\theta,\ueps}$
as in Lemma~\ref{lemma:3I}.
\begin{lemma}\label{lemma:3J}
For each $i=1,\ldots,r$, $\theta\in\Deltad$ and
$\ueps\in\cE_\infty(K)$ the functions $E_\theta$ and
$D_{i,\theta,\ueps}$ are of Iwasawa type.
\end{lemma}
\bPf\ Clearly, $E_\theta(m,\psi)=e_\theta((1+p)^m\psi(1+p)\inv-1)$ where
$e_\theta(T)$ is the power series
$\langle N\ff'(\theta)d_k\rangle(1+T)^{-\ell(N\ff'(\theta)d_k)}\in\bbZ_p[[T]]$.
Thus $E_\theta$ lies in $\IwaZG$. As for $D_{i,\theta,\ueps}$,
if $x\in\bbZ_p^\times$ then, with our identifications, we find
$\talpha_i(x)|_\kinf=\langle x\rangle=\alpha(\langle x\rangle)|_\kinf$.
It follows that the formula $\nu_i(x):=\talpha_i(x)|_\Kinf\alpha(\langle x\rangle)\inv$
defines a homomorphism $\nu_i:\bbZ_p^\times\rightarrow\Delta$
(trivial on $1+p^{n_2}\bbZ_p$, by\refeq{3I}). Thus if
$\chi=\theta\times\psi\in\Ginfd$ then
$(\chi\circ\talpha_i)(x)=
\chi(\talpha_i(x)|_\Kinf)=
(\theta\circ\nu_i)(x)\psi(\langle x\rangle)
$. Therefore, if $(m,\psi)\in\ZGdd$ and $\chi=\theta\times\psi$
we have
$
D_{i,\theta,\ueps}(m,\psi)=
I_{\bbZ_p^\times,\mu'_{i,\theta,\ueps}}(m,\psi)
$
where
$\mu'_{i,\theta,\ueps}\in {\rm Meas}(\bbZ_p^\times, H^{(i)}(\theta))$
denotes the measure
$\langle x\rangle(\theta\circ\nu_i)(x)\mu_{i,\theta,\ueps}(x)$, which depends
only on $i$, $\theta$ and $\ueps$. The result now follows from
Lemma~\ref{lemma:3G}.\ePf
\noindent
To prove Proposition~\ref{prop:3C} we may again assume (by linearity) that
$\ueta$ is as in\refeq{3S.5}. Combining
Equations\refeq{3S} and\refeq{3T} (with $n=n(\psi)$)
and using also\refeq{3N.6} and the definitions
of  $D_{i,\theta,\ueps}$ and $E_\theta$,
we find that for all $\theta\in\Deltad$,
and all $(m,\psi)\in\ZGdd$
\begin{eqnarray*}
\lefteqn{C_{\theta,\ueta}(m,\psi)=
{\displaystyle\frac{2^r}{\sqrt{d_k}}}
\left[
g(\hchi)\inv
\prod_{i=1}^{r}g_{n(\psi)}(\chi\circ\talpha_i)\prod_{i=1}^{r}\theta(\trho_i)^{-n(\psi)}
\psi(\langle N\ff'(\theta)d_k\rangle)\inv
\right]\times}&&\\
&&\hspace*{20em} E_\theta(m,\psi)\det(D_{i,\theta,\ueps_l}(m,\psi))_{i,l=1}^r
\end{eqnarray*}
where, as usual, $\chi$ denotes $\theta\times\psi$.
By Lemma~\ref{lemma:3J}, the proof of Proposition~\ref{prop:3C} will be complete
once we have proven that the quantity in square brackets is of Iwasawa
type, when considered as a function of $(m,\psi)$. However it is evidently
independent of $m$ so this comes down to showing
\begin{prop}\label{prop:3D}
For each fixed
$\theta\in\Deltad$ the quantity
\[
g(\hchi)\inv\prod_{i=1}^{r}g_{n(\psi)}
(\chi\circ\talpha_i)\prod_{i=1}^{r}\theta(\trho_i)^{-n(\psi)}
\psi(\langle N\ff'(\theta)d_k\rangle)\inv
\]
(with $\chi=\theta\times\psi$) is independent of the character
$\psi\in\Gammadd$
\end{prop}
To prove this, we fix $\theta\in\Deltad$,
and suppose that $\chi=\theta\times\psi\in\Ginfdd$ for some $\psi\in\Gammadd$.
Recall that the character $\hchi$ is regarded as
a (primitive) character of the ray-class group $\Cl_{\ff(\chi)}(k)$
with values either in $\bbC$ or (via $j$) in $\bbC_p$
and that the $p$-adic global Gauss sum which we are denoting $g(\hchi)$
is simply the embedding by $j$ of the complex sum attached to $\hchi$
as defined in~\cite[\S 6.4]{twizas}.
In the notation of the latter we therefore have
\[
g(\hchi)=\sum_{a\in R}
j(e({\rm Tr}_{k/\bbQ}(a)))\hchi([a\ff(\chi)\cD]_{\ff(\chi)})\inv
\]
where $e(x)=\exp(2\pi ix)$ and $R$ is any set
of representatives in $k$ for `$\cT(\ff(\chi),\cD\inv)$', namely
the points of $k/\cD\inv$ whose $\cO$-annihilator is exactly $\ff(\chi)$.
As explained in \emph{ibid.},
for any $a\in R$ the ideal $a\ff(\chi)\cD$ is integral and prime to $\ff(\chi)$.
Furthermore its class $[a\ff(\chi)\cD]_{\ff(\chi)}$ in $\Cl_{\ff(\chi)}(k)$ depends only on
the class $a+\cD\inv\in k/\cD\inv$ as, of course, does $e({\rm Tr}_{k/\bbQ}(a))$.
Let us fix, for each $\theta\in\Deltad$, an element $x_\theta$ of $k$ satisfying
\begin{eqnarray}
\ord_{\fq}(x_\theta)&=&-\ord_{\fq}(\ff'(\theta)\cD)\ \mbox{for all primes
$\fq|\ff'(\theta)$,}
\label{eq:3U}\\
\ord_{\fq}(x_\theta)&\geq&0\ \mbox{for all
$\fq\ndiv\ff'(\theta)$ and}
\label{eq:3V}\\
\ord_{\fp_i}(x_\theta-1)&\geq&n_2\ \mbox{for $i=1,\ldots,r$}
\label{eq:3W}
\end{eqnarray}
These conditions imply that
$\min\{\ord_\fq(p^{-n(\psi)-1}x_\theta),\ord_\fq(\cD\inv)\}=
\ord_\fq(p^{-n(\psi)-1}\ff'(\theta)\inv\cD\inv)$ for every prime
$\fq$ of $\cO$. (In the case $\fq|p$, note that not only $\ff'(\theta)$
but also $\cD$ and $x_\theta$ are prime
to $p$, by Conditions\refeq{3C} and\refeq{3W}).
Using\refeq{3N.75} we obtain
\[
p^{-n(\psi)-1}x_\theta\cO+\cD\inv=p^{-n(\psi)-1}\ff'(\theta)\inv\cD\inv=
\ff(\chi)\inv\cD\inv
\]
in other words $p^{-n(\psi)-1}x_\theta$ generates
$\ff(\chi)\inv\cD\inv/\cD\inv$ as a (free, rank-1) $(\cO/\ff(\chi))$-module.
It follows easily that we may take
$R$ to be $p^{-n(\psi)-1}x_\theta R_{\ff(\chi)}$ where
$R_{\ff(\chi)}$ is any set
of representatives in $\cO$ of $(\cO/\ff(\chi))^\times$.
Thus, dropping `$j$' from the notation once more, we can
rewrite $g(\hchi)$ as
\begin{eqnarray}
g(\hchi)&=&\sum_{b\in R_{\ff(\chi)}}
e({\rm Tr}_{k/\bbQ}(p^{-n(\psi)-1}x_\theta b))
\hchi([x_\theta b\ff'(\theta)\cD]_{\ff(\chi)})\inv\nonumber\\
&=&\hchi([x_\theta \ff'(\theta)\cD]_{\ff(\chi)})\inv
\sum_{b\in R_{\ff(\chi)}}
e(p^{-n(\psi)-1}{\rm Tr}_{k/\bbQ}(x_\theta b))
\hchi([b\cO]_{\ff(\chi)})\inv\label{eq:3X}
\end{eqnarray}
It will be necessary to compare this explicit expression for the global
Gauss sum with the product of the local Gauss sums
$g_{n(\psi)}(\chi\circ\talpha_i)$ \etc\ Such a comparison is hard to
find in the literature. To establish the precise connection
we shall therefore recall some basic facts and notations from the
idelic framework of class field theory.
For any finite abelian extension $M$ of $k$ and each
place $v$ of $k$ we write $M^v$ for the completion $M_w$ of
$M$ at some prime $w$ dividing $v$ and
$\varphi_{v}=\varphi_{v,M/k}$ for the local reciprocity map
from $k_v^\times$ to $\Gal(M^v/k_v)$ identified with the
decomposition subgroup $G(M/k)_v\inv$
of $\Gal(M/k)$ above $v$.
This is independent of the choice of $w$ given $v$ and,
as $v$ varies, the maps $\varphi_v$ give rise to a well-defined global
map $\varphi_{M/k}$
from the id\`ele group ${\rm Id}(k)$ of $k$ to $\Gal(M/k)$, by
sending the id\`ele $(a_v)_v$ to $\prod_{v}\varphi_{v,M/k}(a_v)$.
It is well-known that $\varphi_{M/k}$ is surjective and that
$\ker(\varphi_{M/k})$ contains $k^\times$ (embedded diagonally).
If $M/k$ is abelian but infinite, we shall still
denote by $\varphi_{M/k}$ the (not-necessarily-surjective)
map from $\Id(k)$ to $\Gal(M/k)$ which is
obtained as the limit of the $\varphi_{L/k}$ over finite sub-extensions $L/k$.
Similarly, we shall still write
$\varphi_{v,M/k}:k_v^\times\rightarrow G(M/k)_v\inv\subset\Gal(M/k)$
for the limit of the local reciprocity maps $\varphi_{v,L/k}$.
The \emph{content} $c(\underline{a})$
of an id\`ele $\underline{a}=(a_v)_v\in\Id(k)$
is defined to be the fractional ideal
$c(\underline{a})=\prod_{\fq}\fq^{\ord_{\fq}(a_{\fq})}$ where the product runs over
all prime ideals $\fq$ of $\cO$ (identified with the corresponding finite places).
Now let $\ff$ be an ideal of $\cO$. Every $\underline{a}\in\Id(k)$
can be written as $x\underline{d}$ where $x\in k^\times$ and
$\underline{d}=(d_v)_v\in\Id(k)$
satisfies $d_v\in 1+\ff\cO_{k_v}$ for each $v|\ff$.
We define $\gamma_\ff(\underline{a})$ to be the class $[c(\underline{d})]_\ff$ in
$\Cl_\ff(k)$ of the content of $\underline{d}$.
This gives a well-defined, surjective
homomorphism $\gamma_\ff:\Id(k)\rightarrow\Cl_\ff(k)$ for which one checks that the
composite with the Artin isomorphism $\Cl_\ff(k)\rightarrow\Gal(k(\ff)/k)$ is nothing but
$\varphi_{k(\ff)/k}$.

Now the character $\chi=\theta\times\psi\in\Ginfdd$ as above may be inflated
to $\tGinf$ or `descended' to $\Gal(K_\chi/k)$ where $K_\chi\subset k(\ff(\chi))\cap\tKinf$
is the subfield that it cuts out. The above remarks (together
with the compatibilities of $\varphi_{K_\chi/k}$,
$\varphi_{k(\ff(\chi))/k}$,  $\varphi_{\tKinf/k}$ and $\varphi_{\Kinf/k}$
with the various restriction maps) show that, in our notation,
\beql{eq:3extra}
\hchi\circ\gamma_{\ff(\chi)}=
\chi\circ\varphi_{\tKinf/k}=
\chi\circ\varphi_{\Kinf/k}
\eeq
For each $i=1,\ldots,r$ the embedding $\tau_i$ (really
$j\circ\tau_i$) of $k$ in $\bbQ_p$
extends to an isomorphism $k_{\fp_i}\rightarrow\bbQ_p$ which we shall denote
$\hat{\tau}_i$.
\begin{lemma}\label{lemma:3K}
Supppose that $\underline{a}=(a_v)_v\in\Id(k)$ is an id\`ele satisfying
$a_v>0$ for every real place $v$ of $k$. Then $\varphi_{\tkinf/k}(\underline{a})$
acts on $\mu_{p^\infty}$ as
$Nc(\underline{a})\prod_{i=1}^r\hat{\tau}_i(a_{\fp_i})\inv\in\bbZ_p^\times$.
\end{lemma}
\bPf\ We write $\tilde{\bbQ}_\infty$ instead of $\bbQ(\mu_{p^\infty})$.
The properties of $\varphi$ imply that
$\varphi_{\tkinf/k}(\underline{a})|_{\tilde{\bbQ}_\infty}=
\varphi_{\tilde{\bbQ}_\infty/\bbQ}(N_{k/\bbQ}\underline{a})=
\varphi_{\tilde{\bbQ}_\infty/\bbQ}(\underline{s})$ where $\underline{s}=(s_v)_v$
is the rational id\`ele $Nc(\underline{a})\inv N_{k/\bbQ}\underline{a}\in\Id(\bbQ)$. But
$s_q\in\bbZ_q^\times$ for every prime number $q$ and $s_\infty>0$. Therefore
$\varphi_{\tilde{\bbQ}_\infty/\bbQ}(\underline{s})$ acts on $\mu_{p^\infty}$
by $s_p\inv=Nc(\underline{a})\prod_{i=1}^r\hat{\tau}_i(a_{\fp_i})\inv$.
(For the last statement, see \S5.7, Ch.~VII
of~\cite{C-F}.
In fact since every $q\neq p$ is unramified in $\tilde{\bbQ}_\infty/\bbQ$
it follows that $\varphi_{v,\tilde{\bbQ}_\infty/\bbQ}(s_v)=1$ for
each place $v\neq p$. So $\varphi_{\tilde{\bbQ}_\infty/\bbQ}(\underline{s})=
\varphi_{p,\tilde{\bbQ}_\infty/\bbQ}(s_p)$ which indeed acts
as $s_p\inv$ by local class field theory. See
the Remark of \S7,~Ch.~XIV of~\cite{SerreLF} or
\S3.7,~Ch.~VI of~\cite{C-F}.)\ePf
\noindent For each $b\in k^\times$ and $i=1,\ldots,r$ we write $\ubi\in\Id(k)$
for the id\`ele with entry $b\inv$ (\emph{not} $b$)
at the place $\fp_i$ and $1$ everywhere else.
If $b$ is a local unit at $\fp_i$ (\resp\ $b=p$) then the above lemma
shows that $\varphi_{\tkinf/k}(\ubi)$ acts by $\tau_i(b)$ (\resp\ acts trivially) on
$\mu_{p^\infty}$. Moreover, $\varphi_{L^{(i)}/k}(\ubi)=\varphi_{\fp_i,L^{(i)}/k}(b\inv)$
equals $1$ (\resp\ equals the inverse of the Frobenius above $\fp_i$)
by local class field theory, since $L^{(i)}/k$ is unramified at $\fp_i$. We conclude:
\beql{eq:3X.5}\mbox{
$\varphi_{\tKinf/k}(\ubi)=
\talpha_i(\tau_i(b))\in\tGamma{i}$ for all $b$ prime to $\fp_i$, and
$\varphi_{\tKinf/k}(\underline{p}^{(i)})=\trho_i\inv\in\tDelta$}
\eeq
For any $b\in k^\times$ we also define $\tub$ (\resp\ $\ub'$)
to be the id\`ele with entry $b\inv$ at
every place dividing $p\ff'(\theta)$ (\resp\ dividing $\ff'(\theta)$) and
$1$ everywhere else,
so that $\tub=\ub'\prod_{i=1}^r\ubi$.
If $b$ is prime to $\ff(\chi)$ for $\chi\in\Ginfdd$ (\ie\ to $p\ff'(\theta)$)
then it is easy to check that $\gamma_{\ff(\chi)}(\tub)=[b\cO]_{\ff(\chi)}$
so, by\refeq{3extra},
\beql{eq:3Y}
b\in R_{\ff(\chi)}\Rightarrow\hchi([b\cO]_{\ff(\chi)})=\chi(\varphi_{\tKinf/k}(\tub))=
\chi(\varphi_{\tKinf/k}(\ub'))
\prod_{i=1}^r\chi(\varphi_{\tKinf/k}(\ubi))
\eeq
On the other hand, if $b$ is prime to $\ff'(\theta)$ then for every prime ideal
$\fq|\ff'(\theta)$, local class field theory gives $\varphi_{\fq,\Kinf/k}(b\inv)\in
G(\Kinf/k)_\fq^0\subset\Delta$ and, moreover,
if $\ord_\fq(b-1)\geq\ord_\fq(\ff'(\theta))$ then $\varphi_{\fq,\Kinf/k}(b\inv)$
lies in
$G(\Kinf/k)_\fq^{\ord_\fq(\ff'(\theta))}$ hence in $\ker \theta$ by\refeq{3N.5}.
It follows easily that there is a well-defined homomorphism
$\theta':(\cO/\ff'(\theta))^\times\rightarrow\bbC_p^\times$, depending only
on $\theta$, which sends the class $\bar{b}\in(\cO/\ff'(\theta))^\times$
to $\theta(\varphi_{\Kinf/k}(\ub'))=\chi(\varphi_{\Kinf/k}(\ub'))=
\chi(\varphi_{\tKinf/k}(\ub'))$.
Combining\refeq{3X.5} with\refeq{3Y}, we obtain
\beql{eq:3Z}
b\in R_{\ff(\chi)}\Rightarrow\hchi([b\cO]_{\ff(\chi)})=
\theta'(\bar{b})\prod_{i=1}^r(\chi\circ\talpha_i)(\tau_i(b))
\eeq
Note also
that $\varphi_{v,\Kinf/k}(p\inv)$ lies in
$G(\Kinf/k)^0_v\subset \ker(\theta)$
for every $v\ndiv p\ff'(\theta)$ hence so also must
$\varphi_{\Kinf/k}(\underline{\tilde{p}})$. Regarding $\tDelta$
as a subgroup of $\Delta$ by restriction as usual, it follows that
\beql{eq:3A1}
\theta'(\bar{p})=
\theta(\varphi_{\Kinf/k}(\underline{p}'))
=\prod_{i=1}^r \theta(\varphi_{\tKinf/k}(\underline{p}^{(i)}))\inv
=\prod_{i=1}^r \theta(\trho_i)
\eeq
by \refeq{3X.5}.

We turn attention to the first factor of the summand in\refeq{3X}.
Define $f'(\theta)\in\bbZ_{\geq 1}$ by $\ff'(\theta)\cap\bbZ=f'(\theta)\bbZ$.
Since $\ff'(\theta)x_\theta\subset\cD\inv$ (by\refeq{3U} and\refeq{3V})
we may also define a function
\displaymapdef{\xi'=\xi'_{x_\theta}}{\cO/\ff'(\theta)}{\mu_{f'(\theta)}}{\bar{b}}
{j(e({\rm Tr}_{k/\bbQ}(x_\theta b)))}
and set
\[
g(\theta';x_\theta):=\sum_{b\in R_{\ff'(\theta)}}
\xi'(\bar{b})\theta'(\bar{b})\inv
\]
where $R_{\ff'(\theta)}$ is any set of representatives in $\cO$ of
$(\cO/\ff'(\theta))^\times$. (Clearly, $g(\theta';x_\theta)$ depends only on
$\theta$ and the choice of $x_\theta$).
\begin{lemma} For any $b\in\cO$ and $n\in\bbN$, we have
\beql{eq:3B1}
e(p^{-n-1}{\rm Tr}(x_\theta b))=
\xi'(\bar{p}^{-n-1}\bar{b})\prod_{i=1}^r\zeta_n^{\tau_i(x_\theta b)}
\eeq
Where ${\rm Tr}={\rm Tr}_{k/\bbQ}$, $\bar{b}$ denotes
the class of $b$ in $\cO/\ff'(\theta)$, `$e$' stands for
$j\circ e$ and `$\tau_i$' for $j\circ\tau_i$ as usual.
\end{lemma}
\bPf\ Note that $\tau_i(x_\theta b)$ lies in $\bbZ_p$ by\refeq{3V}
so $\zeta_n^{\tau_i(x_\theta b)}$ makes sense. We need to prove
that the two roots of unity, $e(p^{-n-1}{\rm Tr}(x_\theta b))$ and
$\xi'(\bar{p}^{-n-1}\bar{b})\zeta_n^{{\rm Tr}(x_\theta b)}$ are equal.
But their $p^{n+1}$th powers are clearly equal, and so are their $f'(\theta)$th powers
(note that $f'(\theta){\rm Tr}(x_\theta b)$ is an integer). Since
$(p^{n+1},f'(\theta))=1$, the result follows.\ePf
\noindent Now, if $\chi=\theta\times\psi\in\Ginfdd$ then, as
$b$ runs through $R_{\ff(\chi)}$, so the image of $b$
in $(\cO/\ff'(\theta))^\times$ and the images of $\tau_i(b)$
(for $r=1,\ldots,r$) in
$\bbZ_p^\times/(1+p^{n(\psi)+1}\bbZ_p)$ all
run through these quotients independently of each other.
Therefore, fixing any
sets $R_{\ff'(\theta)}$ and $R_{n(\psi)}$ of
representatives for these quotients and
substituting\refeq{3Z} and\refeq{3B1} (with $n=n(\psi)$)
into the R.H.S. of\refeq{3X} we obtain
\begin{eqnarray}
g(\hchi)&=&\hchi([x_\theta \ff'(\theta)\cD]_{\ff(\chi)})\inv
\sum_{b\in R_{\ff(\chi)}}
\left[
\xi'(\bar{p}^{-n(\psi)-1}\bar{b})\theta'(\bar{b})\inv
\prod_{i=1}^r\left(\zeta_{n(\psi)}^{\tau_i(x_\theta b)}
(\chi\circ\talpha_i)(\tau_i(b))\right)\inv
\right]\nonumber\\
&=&\hchi([x_\theta \ff'(\theta)\cD]_{\ff(\chi)})\inv
\left(\sum_{b\in R_{\ff'(\theta)}}
\xi'(\bar{p}^{-n(\psi)-1}\bar{b})\theta'(\bar{b})\inv\right)
\prod_{i=1}^r
\sum_{u\in R_{n(\psi)}}\zeta_{n(\psi)}^{\tau_i(x_\theta)u}
(\chi\circ\talpha_i)(u)\inv\nonumber\\
&=&\hchi([x_\theta \ff'(\theta)\cD]_{\ff(\chi)})\inv
\theta'(\bar{p})^{-n(\psi)-1}
g(\theta';x_\theta)
\prod_{i=1}^r[(\chi\circ\talpha_i)
(\tau_i(x_\theta))g_{n(\psi)}(\chi\circ\talpha_i)]\label{eq:3C1}
\end{eqnarray}
(which shows incidentally that $g(\theta';x_\theta)\neq 0$).
Rearranging\refeq{3C1} and using\refeq{3A1} and the
fact (from \refeq{3I} and\refeq{3W})
that $(\chi\circ\talpha_i)(\tau_i(x_\theta))=
\psi(\tau_i(x_\theta))$ for each $i$, we find
\begin{eqnarray}
\lefteqn{g(\hchi)\inv\prod_{i=1}^{r}g_{n(\psi)}
(\chi\circ\talpha_i)\prod_{i=1}^{r}\theta(\trho_i)^{-n(\psi)}
\psi(\langle N\ff'(\theta)d_k\rangle)\inv=}&&\hspace{20em}\nonumber\\
&&
\left(\prod_{i=1}^{r}\theta(\trho_i)\right)
g(\theta';x_\theta)\inv
\left[\hchi([x_\theta \ff'(\theta)\cD]_{\ff(\chi)})
\psi(N_{k/\bbQ}x_\theta)\inv
\psi(\langle N\ff'(\theta)d_k\rangle)\inv\right]\label{eq:3XS}
\end{eqnarray}
Finally, for each $\theta\in\Deltad$ we fix an id\`ele
$\underline{y}_\theta=(y_{\theta,v})_v\in\Id(k)$
such that $y_{\theta,v}=1$ for
all places $v$ which are infinite or divide $p\ff'(\theta)$
and $c(\underline{y}_\theta)=x_\theta\ff'(\theta)\cD$.
Then
\beql{eq:3D1}
\hchi([x_\theta \ff'(\theta)\cD]_{\ff(\chi)})=
\hchi(\gamma_{\ff(\chi)}(\underline{y}_\theta))=
\chi(\varphi_{\Kinf/k}(\underline{y}_\theta))
\eeq
by\refeq{3extra}.
Also, keeping in mind our identification of
$\tiGamma$ with $\bbZ_p^\times$ via $\kappa$
\etc, Lemma~\ref{lemma:3K} gives
\[
\varphi_{\Kinf/k}(\underline{y}_\theta)|_\kinf=
\langle\kappa(\varphi_{\tkinf/k}(\underline{y}_\theta))\rangle=
\langle|N_{k/\bbQ}x_\theta|N\ff'(\theta)d_k\rangle=
N_{k/\bbQ}x_\theta\langle N\ff'(\theta)d_k\rangle\in\Gamma
\]
so that the element $z_\theta:=
\varphi_{\Kinf/k}(\underline{y}_\theta)
\alpha(N_{k/\bbQ}x_\theta\langle N\ff'(\theta)d_k\rangle)\inv$ lies in
$\Delta$. Using Equation\refeq{3D1}, it follows that the
third factor on the R.H.S. of\refeq{3XS} equals
$\chi(z_\theta)=\theta(z_\theta)$ and so it depends only
on $\theta$, not on $\psi$ (in fact,
it is even independent of the choice of $\underline{y}_\theta$
given $\theta$).
But the other two factors on the R.H.S. of\refeq{3XS}
are also independent of $\psi$
so we have proven Proposition~\ref{prop:3D},
hence also Proposition~\ref{prop:3C}, hence also Theorem~\ref{thm:3A}.
\section{Semilocal Weakenings of $P1$ and $P2$}\label{section:semiloc}
The aim of this section is to explain
the weaker `semilocal' variants  of properties $P1$ and $P2$
namely $\widehat{P1}$ and $\widehat{P2}$ which were
mentioned in the introduction.
We shall also indicate how the constructions and results of
Section~\ref{section:consP2} may be adapted
to obtain (for example) an analogue of Corollary~\ref{cor:3A}
which assumes only $\widehat{P2}$.

For fixed $p\neq 2$ and any number field $F$, we write $U(F)$ for the
$p$-semilocal unit group $(\bbZ_p\otimes_\bbZ\cO_F)^\times$
and $U^{(1)}(F)$ for its Sylow pro-$p$ subgroup. We identify these
respectively with the products $\prod_{\fP|p}U(F_\fP)$ and
$\prod_{\fP|p}U^{(1)}(F_\fP)$ of the local (\resp\ the principal
local) units in the completions of $F$ at all primes $\fP$ dividing $p$.
Any embedding $\iota:F\rightarrow\bbC_p$
defines a natural homomorphism
$U(F)\rightarrow U(\bbC_p)$ which we shall denote
$\widehat{\iota}$. (In the product realisation, $\iota$ extends to an
embedding of $F_\fP\rightarrow\bbC_p$
for a unique prime $\fP|p$ and $\widehat{\iota}$ is
obtained by composing this map with the projection
$U(F)\rightarrow U(F_\fP)$.)
Also, we shall denote by $\cU_\infty^{(1)}(F)$ the inverse limit of $U^{(1)}(F_n)$
over $n\geq n_1(F,p)$ with respect to the maps
$U^{(1)}(F_n)\rightarrow U^{(1)}(F_m)$ coming from the products of the
local norms.

Now suppose $K/k$ and $p$
satisfy\refeq{2A},\refeq{2B} and\refeq{2C} then, for each
$n\geq n_1$ we endow $U^{(1)}(K_n)$ with the
usual structure of $\bbZ_p G_n$-module
and let $a_n:E(K_n)\rightarrow U^{(1)}(K_n)$
denote the $\bbZ G_n$-homomorphism which is
the composite of the natural map $E(K_n)\rightarrow U(K_n)$ with the usual
projection of $U(K_n)$ on $U^{(1)}(K_n)$. (Thus $\ker(a_n)=\{\pm 1\}$ and
Leopoldt's conjecture -- which we need not assume here --
predicts that ${\rm Im}(a_n)$ spans a $\bbZ_p$-submodule
of co-rank $1$ in $U^{(1)}(K_n)$.)
The map $\bigwedge^r a_n$ gives rise to an
obvious $\bbQ G_n$-linear map $b_n$, say, from
$\bigwedge^r_{\bbQ G_n}\bbQ E(K_n)=\bbQ\otimes_\bbZ\bigwedge^r_{\bbZ G_n} E(K_n)$
to $\bbQ_p \otimes_{\bbZ_p}\bigwedge^r_{\bbZ_p G_n} U^{(1)}(K_n)$.
Now $\cU_\infty^{(1)}(K)$
is naturally a module for
the \emph{completed group-ring} $\bbZ_p[[G_\infty]]$ (namely, the
inverse limit over $n\geq n_1$ of the rings $\bbZ_p G_n$ w.r.t.\ the natural restriction
maps $\bbZ_p G_n\rightarrow \bbZ_p G_m$). We may therefore form the exterior
power $\bigwedge^r_{\bbZ_p[[G_\infty]]}\cU_\infty^{(1)}(K)$ which is equipped
with a natural map $\hat{\beta}_n$ to $\bigwedge^r_{\bbZ_p G_n} U^{(1)}(K_n)$
for each $n\geq n_1$. Furthermore,
the map $a_\infty={\displaystyle\lim_{\leftarrow}}\,a_n:
\cE_\infty(K)\rightarrow\cU_\infty^{(1)}(K)$ naturally defines
a homomorphism $b_\infty$ from $\bigwedge^r_{\bbZ G_\infty}\cE_\infty(K)$
to $\bigwedge^r_{\bbZ_p[[G_\infty]]}\cU_\infty^{(1)}(K)$.
We now have the commutative diagram~(\ref{diag:2})
\begin{figure}
\setlength{\unitlength}{5.5mm}
\beql{diag:2}
\begin{picture}(22,10)(0,0)
\put(8.5,1){\vector(1,0){5.3}}
\put(8.9,9){\vector(1,0){4.9}}
\put(1.7,8.4){\vector(0,-1){6.8}}
\put(18,8.2){\vector(0,-1){6.5}}
\put(2,5){\makebox(0,0)[l]{$ \beta_n$}}
\put(18.3,5){\makebox(0,0)[l]{$\hat{\beta}_n$}}
\put(11,1.2){\makebox(0,0)[b]{$ b_n $}}
\put(11,9.2){\makebox(0,0)[b]{$b_\infty$}}
\put(2.5,1){\makebox(0,0){$\bigwedge^r_{\bbQ G_n}\bbQ E(K_n)=
     \bbQ\otimes_\bbZ\bigwedge^r_{\bbZ G_n} E(K_n)$}}
\put(2.5,9){\makebox(0,0){$\bigwedge^r_{\bbQ G_\infty}\bbQ\cE_\infty(K)=
          \bbQ\otimes_\bbZ\bigwedge^r_{\bbZ G_\infty}\cE_\infty(K)$}}
\put(18,1){\makebox(0,0){$\bbQ_p\otimes_{\bbZ_p}\bigwedge^r_{\bbZ_pG_n}U^{(1)}(K_n)$}}
\put(18,9){\makebox(0,0){
           $\bbQ_p\otimes_{\bbZ_p}\bigwedge^r_{\bbZ_p[[G_\infty]]}\cU_\infty^{(1)}(K)$}}
\end{picture}
\eeq
\end{figure}
where we have extended
$\hat{\beta}_n\ \forall n\geq n_1$ and $b_\infty$
linearly to the tensor products. We formulate
\begin{property}{$\widehat{P1}(K/k,p)$}
There exists an element
$\hat{\ueta}\in\bbQ_p\otimes_{\bbZ_p}
\bigwedge^r_{\bbZ_p[[G_\infty]]}\cU_\infty^{(1)}(K)$
such that, for every $n\geq n_1$
we have
$\hat{\beta}_n(\hat{\ueta})=b_n(\eta_n)$
where  $\eta_n$ lies in $V_n^0$
and is a solution of
$C1(K_n/k,\emptyset)$ (the unique solution in $V_n^0$) \ie\
\[
\Phi_n(1)=\frac{2^r}{\sqrt{d_k}}
R_n(\eta_n)\ \ \ \forall n\geq n_1
\]
(We say that $\hat{\ueta}$
\emph{demonstrates} $\widehat{P1}(K/k,p)$).
\end{property}
\noindent\rem\label{rem:P1vsP1hat} This property is implied by $P1(K/k,p)$. Indeed, if
$\ueta$ demonstrates $P1(K/k,p)$ then clearly $b_\infty(\ueta)$
demonstrates $\widehat{P1}(K/k,p)$. However the converse implication
is not obviously true. For example, $\widehat{P1}(K/k,p)$
at most implies that the \emph{$p$-part}
of the denominator of $\eta_n$ is bounded as $b\rightarrow\infty$.
(Compare Remark~\ref{rem:denoms}).\vspace{1ex}\\
One way to formulate $\widehat{P2}(K/k,p)$
would simply be to replace the above
condition that $\eta_n\in V_n^0$ be a solution of $C1(K_n/k,\emptyset)\ \forall n\geq n_1$in
by the requirement that $p^re_n\eta_n$ be a solution of $C2(K_n/k,T_p,p)\ \forall n\geq n_2$,
or indeed the
equivalent equation involving $R_{n,p}(e_n\eta_n)$. There is, however, a neater, equivalent
formulation, using the fact that the map $R_{n,p}$ can be `extended' to
$\bbQ_p\otimes_{\bbZ_p}\bigwedge^r_{\bbZ_p[G_n]}U^{(1)}(K_n)$
(\ie\ it factors through $b_n$). Indeed, for each $i=1,\ldots,r$ we may
`extend' $\lambda_{K_n/k,i,p}$ to a $\bbZ_p G_n$-linear  map
$\hat{\lambda}_{K_n/k,i,p}:U^{(1)}(K_n)\rightarrow\bbC_p$ by replacing the embedding
$j\ttau_i$ with $\widehat{j\ttau_i}$ in the definition of $\lambda_{K_n/k,i,p}$.
Then we obtain a unique
$\bbQ_p G_n$-linear
regulator map $\hat{R}_{n,p}=\hat{R}_{K_n/k,p}$ from
$\bbQ_p \otimes_{\bbZ_p}\bigwedge^r_{\bbZ_p G_n} U^{(1)}(K_n)$ to $\bbC_p G_n$
by setting
$\hat{R}_{n,p}(1\otimes (u_1\wedge\ldots\wedge u_r)):=
\det(\hat{\lambda}_{K_n/k,i,p}(u_s))_{i,s=1}^r$.
It should be clear that $\hat{R}_{n,p}\circ b_n=R_{n,p}$.
For all $K/k$ and $p$
satisfying\refeq{2A},\refeq{2B} and\refeq{2C} we formulate
\begin{property}{$\widehat{P2}(K/k,p)$}
There exists an element $\hat{\ueta}\in\bbQ_p\otimes_{\bbZ_p}
\bigwedge^r_{\bbZ_p[[G_\infty]]}\cU_\infty^{(1)}(K)$
such that, for every $n\geq n_2$,
$\hat{\beta}_n(\hat{\ueta})$ lies in $b_n(V_n^0)$ and satisfies
\[
\Phi_{n,p}(1)=\frac{2^r}{j(\sqrt{d_k})}
\hat{R}_{n,p}(e_n\hat{\beta}_n(\hat{\ueta}))
\]
(We say that $\hat{\ueta}$
\emph{demonstrates} $\widehat{P2}(K/k,p)$).
\end{property}
\rem\ The exact analogue of Proposition~\ref{prop:2B} now holds, namely
one may replace $P1$ by $\widehat{P1}$ and $P2$ by $\widehat{P2}$
(and $\ueta$ by $\hat{\ueta}$) in the statement of this proposition,
the proof remaining essentially the same.\vspace{1ex}\\
\noindent Next, we strengthen the hypothesis\refeq{2C}
to\refeq{3C} and sketch the analogues  of some of the constructions
of Section~\ref{section:consP2}. Just as before,
for any $i=1,\ldots,r$ and  $\underline{u}=(u_n)_n\in\cU_\infty^{(1)}(K)$,
we use the norm maps to $\tL{i}_{n}$ for $N>n\geq 0$
to extend the subsequence $(u_n)_{n\geq N}$
(where $N=\max(n_1,n_2-1)$) to an element of
$\cU_\infty^{(1)}(\tL{i}_0)$. Then we
apply $\widehat{j\ttau_i}$ to get an element
of $U_\infty(H^{(i)})$. The analogue of Lemma~\ref{lemma:3C}
obviously holds with $\underline{u}$ and $\widehat{j\ttau_i}$
in place of  $\underline{v}$ and $j\ttau_i$.
Thus for any $(t,n)\in\ZZnt$ we get a well-defined homomorphism
$\hat{\cL}_{i,t,n}:\cU_\infty^{(1)}(K)\rightarrow\cO_{\tH^{(i)}_n}G_n$
by simply replacing
$j\ttau_i$ by $\widehat{j\ttau_i}$ (and $\ueps\in\cE_\infty(K)$ by
$\underline{u}\in\cU_\infty^{(1)}(K)$) in the definition
of $\check{\cL}_{i,t,n}$. The analogue of Equation\refeq{3F} holds with
$\check{\cL}_{i,t,n}$, $\lambda_{K_n/k,i,p}$ and $\ueps$ replaced
by $\hat{\cL}_{i,t,n}$, $\hat{\lambda}_{K_n/k,i,p}$ and
$\uu$ respectively. Moreover, we have
\begin{lemma}
For each $i,n,t$ as above, the map $\hat{\cL}_{i,t,n}$ is
$\bbZ_p[[G_\infty]]$-linear with
$\bbZ_p[[\Ginf]]$ acting on
the group-ring $\cO_{\tH^{(i)}_n}G_n$ through the quotient $\bbZ_pG_n$.
\end{lemma}
\bPf\ (Sketch.) $\hat{\cL}_{i,t,n}$ is clearly $\bbZ_p$-linear and also
$G_\infty$-linear with $G_\infty$ acting through the $G_n$. In particular,
it is $\bbZ_p G_\infty$-linear. Let $I_n=\ker(\bbZ_p[[G_\infty]]\rightarrow\bbZ_pG_n)$.
Then $\bbZ_p[[G_\infty]]=\bbZ_p G_\infty+I_n$ so it suffices to
prove that $\hat{\cL}_{i,t,n}(z\uu)=0$ for all $\uu\in\cU_\infty^{(1)}(K),\ z\in I_n$.
But $n\geq n_2$ implies that
$\Gal(\Kinf/K_n)$ is isomorphic to $\bbZ_p$, with
topological generator $\gamma_n$, say, and a
compactness argument (for example) shows
that $I_n=(\gamma_n-1)\bbZ_p[[\Ginf]]$. Thus
$\hat{\cL}_{i,t,n}(z\uu)=\hat{\cL}_{i,t,n}((\gamma_n-1)z'\uu)=
\hat{\cL}_{i,t,n}(\gamma_nz'\uu)-\hat{\cL}_{i,t,n}(z'\uu)=0$ since the
image of $\gamma_n$ in $G_n$ is $1$.\ePf
It follows that for each $(t,n)\in\ZZnt$ there is a unique, well-defined
$\bbQ_p\otimes_{\bbZ_p}\bbZ_p[[\Ginf]]$-linear
higher regulator map $\hat{\cR}_{t,n}$ from
$\bbQ_p \otimes_{\bbZ_p}\bigwedge^r_{\bbZ_p[[\Ginf]]}\cU_\infty^{(1)}(K)$ to
$\bbC_p G_n$ satisfying
$\hat{\cR}_{t,n}(1\otimes (\uu_1\wedge\ldots\wedge \uu_r)):=
\det(\hat{\cL}_{i,t,n}(\uu_s))_{i,s=1}^r$.
It should be clear that $\hat{\cR}_{t,n}$ `extends'
$\check{\cR}_{t,n}$ in the sense that
$\hat{\cR}_{t,n}\circ b_\infty=\check{\cR}_{t,n}$. Moreover, the analogue of
Lemma~\ref{lemma:3D} holds, namely
\beql{eq:3E1}
\hat{\cR}_{0,n}(\hat{\ueta})=\hat{R}_{n,p}(e_n\hat{\beta}_n(\hat{\ueta}))\ \ \
\forall\,\hat{\ueta}\in\bbQ_p \otimes_{\bbZ_p}
{\textstyle\bigwedge^r_{\bbZ_p[[G_\infty]]}} \cU^{(1)}_\infty(K),\ \forall\,n\geq n_2
\eeq
The (hypothetical) relation\refeq{3G}
is of course replaced by the following one, for
$\hat{\ueta}\in\bbQ_p \otimes_{\bbZ_p}
{\textstyle\bigwedge^r_{\bbZ_p[[G_\infty]]}} \cU^{(1)}_\infty(K)$ and
$(m,n)\in\bbZ\times\bbZ_{\geq n_2}$:
\beql{eq:3F1}
\Phi_{n,p}(m)=\frac{2^r}{j(\sqrt{d_k})}\hat{\cR}_{1-m,n}(\hat{\ueta})
\eeq
and the most important thing to check
is that Theorem~\ref{thm:3A} now extends as follows.
\begin{thm}\label{thm:4A}
Suppose that $K$ and $p$ satisfy\refeq{2A},\refeq{2B} and\refeq{3C}.
For a given element $\hat{\ueta}\in\bbQ_p \otimes_{\bbZ_p}
{\textstyle\bigwedge^r_{\bbZ_p[[G_\infty]]}} \cU^{(1)}_\infty(K)$,
Equation\refeq{3F1} holds for all pairs $(m,n)\in\bbZ\times \bbZ_{\geq n_2}$
if and only if it holds for infinitely many such pairs.
\end{thm}
\bPf\ The key Lemma~\ref{lemma:3I} in the proof of Theorem~\ref{thm:3A}
works with $\uu\in\cU_\infty^{(1)}(K)$ and $\hat{\cL}_{i,1-m,n}$
in place of $\ueps\in\cE_\infty(K)$ and $\check{\cL}_{i,1-m,n}$. Indeed,
the construction of the required measure
(now $\mu_{i,\theta,\uu}$) works exactly as before because,
crucially, it depends only on the power series attached
norm-coherent sequences of \emph{local} units (now those coming from
$\widehat{j\ttau_i}$ applied to $\uu$ and its various Galois conjugates).
Having thus `extended' this lemma, the other aspects of the
proof Theorem~\ref{thm:3A} go across identically, \emph{mutatis mutandis},
as the reader may verify.\ePf
\noindent Finally, from Equation\refeq{3E1} and
Theorem~\ref{thm:4A} we deduce the following analogue of Corollary~\ref{cor:3A}
\begin{cor}\label{cor:4A}
If~\refeq{2A},\refeq{2B} and\refeq{3C} hold and
$\hat{\ueta}\in\bbQ_p \otimes_{\bbZ_p}
{\textstyle\bigwedge^r_{\bbZ_p[[G_\infty]]}} \cU^{(1)}_\infty(K))$
demonstrates $\widehat{P2}(K/k,p)$ then it satisfies Equation\refeq{3F1}
for all $(m,n)\in\bbZ\times \bbZ_{\geq n_2}$.\ePf
\end{cor}

\end{document}